\documentclass[12pt]{amsart}

\usepackage{amsfonts}
\usepackage{amssymb}
\usepackage{latexsym}
\usepackage{graphics}
\usepackage[2emode]{psfrag}
\usepackage{mathrsfs}
\usepackage{amsthm}
\usepackage{amsmath}
\usepackage[all]{xy}
\usepackage{enumerate}
\usepackage{hyperref}
\usepackage{stmaryrd}

\DeclareMathOperator{\newsquare}{\raisebox{0ex}{\scalebox{0.75}{\ensuremath{\square}}}}
\DeclareMathOperator{\newblacksquare}{\raisebox{0ex}{\scalebox{0.75}{\ensuremath{\blacksquare}}}}

\DeclareMathOperator{\newtriangle}{\raisebox{0ex}{\scalebox{0.75}{\ensuremath{\triangledown}}}}
\DeclareMathOperator{\newblacktriangle}{\raisebox{0ex}{\scalebox{0.75}{\ensuremath{\blacktriangledown}}}}
\DeclareMathOperator{\newtriangled}{\raisebox{0.1ex}{\scalebox{0.5}{\ensuremath{\triangledown}}}}
\DeclareMathOperator{\newblacktriangled}{\raisebox{0.1ex}{\scalebox{0.5}{\ensuremath{\blacktriangledown}}}}
\DeclareMathOperator{\newdiamond}{\raisebox{0.2ex}{\scalebox{0.75}{\ensuremath{\lozenge}}}}
\DeclareMathOperator{\newblackdiamond}{\raisebox{0.2ex}{\scalebox{0.75}{\ensuremath{\blacklozenge}}}}
\DeclareMathOperator{\newdiamondd}{\raisebox{0.2ex}{\scalebox{0.5}{\ensuremath{\lozenge}}}}
\DeclareMathOperator{\newblackdiamondd}{\raisebox{0.2ex}{\scalebox{0.5}{\ensuremath{\blacklozenge}}}}
\DeclareMathOperator{\newtriangleup}{\raisebox{0ex}{\scalebox{0.75}{\ensuremath{\triangle}}}}
\DeclareMathOperator{\newblacktriangleup}{\raisebox{0ex}{\scalebox{1}{\ensuremath{\blacktriangle}}}}

\def\equal#1{\smash{\mathop{=}\limits^{#1}}}

\newcommand{\thlabel}[1]{\label{th:#1}}
\newcommand{\thref}[1]{Theorem~\ref{th:#1}}
\newcommand{\selabel}[1]{\label{se:#1}}
\newcommand{\seref}[1]{Section~\ref{se:#1}}
\newcommand{\lelabel}[1]{\label{le:#1}}
\newcommand{\leref}[1]{Lemma~\ref{le:#1}}
\newcommand{\prlabel}[1]{\label{pr:#1}}
\newcommand{\prref}[1]{Proposition~\ref{pr:#1}}
\newcommand{\colabel}[1]{\label{co:#1}}
\newcommand{\coref}[1]{Corollary~\ref{co:#1}}
\newcommand{\relabel}[1]{\label{re:#1}}
\newcommand{\reref}[1]{Remark~\ref{re:#1}}

\newcommand{\delabel}[1]{\label{de:#1}}
\newcommand{\deref}[1]{Definition~\ref{de:#1}}
\newcommand{\eqlabel}[1]{\label{eq:#1}}
\newcommand{\equref}[1]{(\ref{eq:#1})}

\newcommand{\Hom}{{\rm Hom}}
\newcommand{\End}{{\rm End}}

\newcommand{\im}{{\rm Im}\,}

\newcommand{\Bim}{{\sf Bim}}
\newcommand{\Frm}{{\sf Frm}}
\newcommand{\CAT}{{\sf CAT}}

\newcommand{\Rat}{{\rm Rat}}

\def\tildej{\tilde{\jmath}}
\def\barj{\bar{\jmath}}

\def\ot{\otimes}
\def\bul{\bullet}
\def\ubul{\underline{\bullet}}

\def\id{\textrm{{\small 1}\normalsize\!\!1}}

\def\MM{{\mathbb M}}
\def\NN{{\mathbb N}}

\def\YY{{\mathbb Y}}

\def\cc{{\mathfrak C}}

\def\aaa{{\mathfrak a}}
\def\bbb{{\mathfrak b}}
\def\ccc{{\mathfrak c}}

\def\fff{{\mathfrak f}}

\def\mmm{{\mathfrak m}}
\def\nnn{{\mathfrak n}}

\newcommand{\Aa}{\mathcal{A}}
\newcommand{\Bb}{\mathcal{B}}
\newcommand{\Cc}{\mathcal{C}}
\newcommand{\Dd}{\mathcal{D}}
\newcommand{\Ee}{\mathcal{E}}
\newcommand{\Ff}{\mathcal{F}}
\newcommand{\Gg}{\mathcal{G}}
\newcommand{\Hh}{\mathcal{H}}
\newcommand{\Ii}{\mathcal{I}}
\newcommand{\Jj}{\mathcal{J}}
\newcommand{\Kk}{\mathcal{K}}
\newcommand{\Ll}{\mathcal{L}}
\newcommand{\Mm}{\mathcal{M}}
\newcommand{\Nn}{\mathcal{N}}

\newcommand{\Rr}{\mathcal{R}}
\newcommand{\Ss}{\mathcal{S}}

\newcommand{\Vv}{\mathcal{V}}

\newcommand{\Xx}{\mathcal{X}}

\def\*C{{}^*\hspace*{-1pt}{\Cc}}
\def\*c{{}^*\hspace*{-1pt}{\cc}}

\def\text#1{{\rm {\rm #1}}}

\def\ul{\underline}
\def\dul#1{\underline{\underline{#1}}}
\def\Nat{\dul{\rm Nat}}
\def\Set{\dul{\rm Set}}

\newtheorem{proposition}{Proposition}[section] 
\newtheorem{lemma}[proposition]{Lemma}
\newtheorem{corollary}[proposition]{Corollary}
\newtheorem{theorem}[proposition]{Theorem}

\theoremstyle{definition}
\newtheorem{Definition}[proposition]{Definition}
\newtheorem{Definitions}[proposition]{Definitions}
\newtheorem{example}[proposition]{Example}

\theoremstyle{remark}

\newtheorem{remark}[proposition]{Remark}

\begin{document}

\title{Cofrobenius corings and adjoint functors}

\author[M.C. Iovanov]{M.C. Iovanov${}^*$}
\address{Faculty of Mathematics and Informatics, University of Bucharest,
Str. Academiei 14, RO-010014, Bucharest 1, Romania and\\
State University of New York - Buffalo, 244 Mathematics Building, Buffalo NY, 14260-2900, USA}
\email{yovanov@gmail.com}

\author{J. Vercruysse}    
 \address{Faculty of Engineering, Vrije Universiteit Brussel, Pleinlaan 2, B-1050,
  Brussel, Belgium} 
 \email{joost.vercruysse@vub.ac.be}   
 \urladdr{http://homepages.vub.ac.be/\~{}jvercruy/}
 \thanks{Research supported by the bilateral project BWS04/04 ``New techniques in 
Hopf algebras and graded ring theory" of the Flemish and Romanian 
Governments.}
 \thanks{${}^*$ The author was also partially supported by the contract nr. 24/28.09.07 with UEFISCU "Groups, quantum groups, corings and representation theory" of CNCIS (ID\_1002)}

 \date{\today} 

\begin{abstract}
We study co-Frobenius and more generally quasi-co-Frobenius corings over arbitrary base rings and over PF base rings in particular. We generalize some results about co-Frobenius and quasi-co-Frobenius coalgebras to the case of non-commutative base rings and give several new characterizations for co-Frobenius and more generally quasi-co-Frobenius corings, some of them are new even in the coalgebra situation. 
We construct Morita contexts to study Frobenius properties of corings and a second kind of Morita contexts to study adjoint pairs. Comparing both Morita contexts, we obtain our main result that characterizes quasi-co-Frobenius corings in terms of a pair adjoint functors $(F,G)$ such that $(G,F)$ is locally quasi-adjoint in a sense defined in this note. 
\end{abstract}

\maketitle

\section{introduction}

In the theory of Hopf algebras, quantum groups and their (co)representations, a variety of algebraic structures and corresponding (co)representations have been introduced and studied during the last decades. Among these are comodule (co)algebras and module (co)algebras, the category of Yetter Drinfel'd modules, Doi-Koppinen data or more generally, entwining structures. Although corings were initially introduced by Sweedler \cite{Sw}, they haven't been studied thoroughly until the last decade. The renewed interest in corings has started with an observation made by Takeuchi in \cite{Tak} that corings and their comodules generalize these entwining structures and their entwined modules (see also \cite{Brz:strcor}), and much attention has been devoted to the subject eversince. Moreover, comodules over corings not only generalize many structures important for Hopf algebras and quantum group theory, but they also generalize other important structures such as modules over algebras and comodules over coalgebras, graded modules over graded rings and also, perhaps surprisingly, the chain complexes of modules over an arbitrary ring. Thus, corings and their comodules offer a unifying context for all these structures.

Frobenius and co-Frobenius coalgebras and Hopf algebras, Frobenius ring extensions and Frobenius bimodules have been intensively studied over the last decades. As in other instances, corings offer a general framework for the study of all these Frobenius type properties. For example, the characterization of Frobenius (co)algebras, or Frobenius extensions of rings can be obtained from the more general characterization of Frobenius corings (see \cite{Brz:strcor}). Furthermore, in \cite{BGT, BrzKadWis} the close relations between Frobenius extensions, Frobenius bimodules and Frobenius corings is discussed. 

Although the name indicates differently, the co-Frobenius property of a coring (or coalgebra) is a weakening and not a dualization of the Frobenius property. In particular, although the Frobenius property is left-right symmetric, the co-Frobenius property is not. Nevertheless, coalgebras over a base field which are at the same time left and right co-Frobenius can be understood as a dual version of Frobenius algebras. Indeed, a $k$-algebra $A$ is Frobenius (i.e.\ $A\simeq A^*$ as left, or equivalently right $A$-modules) if and only if the functors ${\rm Hom}_A(-,A)$ and ${\rm Hom}_k(-,k)$ from ${\mathcal M}_A$ to ${}_A{\mathcal M}$ are naturally isomorphic (see \cite{CR}). Similarly, for a coalgebra $C$, considering the natural dual comodule $Rat(C^*_{C^*})$ of $C^C$, it has been recently shown in \cite{mio} that $C$ is left and right co-Frobenius if and only if $C\simeq Rat(C^*_{C^*})$ in ${\mathcal M}^C$ and this allows a functorial-categorical interpretation of this concept: $C$ left and right co-Frobenius if and only if the functors ${\rm Hom}_{C^*}(-,C^*)$ and ${\Hom}_{k}(-,k)$ from ${\mathcal M}^C$ to ${\mathcal M}_{C^*}$ are isomorphic. 

Frobenius corings have a very nice characterization in terms of Frobenius functors. This result says that an $A$-coring ${\cc}$ is Frobenius if and only if the forgetful functor $F:{\mathcal M}^{\cc}\rightarrow {\mathcal M}_A$ is at the same time a left and right adjoint for the induction functor $-\otimes_A{\cc}$. An overview of most results regarding this subject can be found in \cite{CMZ}. A similar categorical interpretation for one sided co-Frobenius and one sided quasi-co-Frobenius coalgebras and corings has remained somewhat mysterious and comes under attention within the theory of corings.

The goal of this paper is to provide this categorical description of quasi-co-Frobenius corings; we also generalize some results of \cite{mio}. The main idea and tool for this will be the construction of several Morita contexts and the interpretation of the (quasi-)co-Frobenius properties in terms of these contexts. Starting from the observation (see \cite[Remarks p 389, Examples 1.2]{Mul}) that a Morita context can be identified with a ($k$-linear) category with two objects, we construct a Morita context relating a coring $\cc$ with its dual $\cc^*$. This context describes the Frobenius property of the coring. We show that if there exists a pair of invertible elements in this Morita context, then the coring is exactly a Frobenius coring. A similar Morita context relates representable functors, such as those used in \cite{mio} and \cite{CR} to describe (co-)Frobenius properties.
A last type of Morita contexts, that is constructed in a different way, describes the adjunction property of a pair of functors. More precise, if there exists a pair of invertible elements in this Morita context, then the pair of functors is exactly an adjoint pair. 
By relating these Morita contexts with (iso)morphisms of Morita contexts, we recover the result that a coring is Frobenius if and only if the forgetful functor and the induction functor make up a Frobenius pair if and only if certain representable functors are isomorphic (\coref{Frobenius}). In particular, using these general Morita contexts, we can formulate a categorical interpretation of co-Frobenius corings and more generally quasi-co-Frobenius corings (see \thref{QcF}). 

The advantage of our presentation is that it clarifies underlying relations between the different equivalent descriptions of the Frobenius property of a coring. In particular, we can explain why the Frobenius property is left-right symmetric and the quasi-co-Frobenius property is not: this is due to a symmetry between several Morita contexts we construct (surfacing in our theory as an isomorphism of Morita contexts) and this symmetry breaks down on the `quasi'-level (see \reref{symmetry}).
A second benefit of our approach is that these Morita contexts and their interrelationship, do exist even if the coring is not Frobenius. This allows us to give a categorical interpretation of co-Frobenius and quasi-co-Frobenius corings by means of these Morita contexts.

This paper is organized as follows. We recall some preliminary results about corings and comodules in \seref{preliminaries}. In \seref{locproj} we give a new characterization for locally projective modules in the sense of Zimmermann-Huisgen \cite{ZH} and rings with local units. 
In \seref{adjoint} we give a new interpretation to the notion of adjoint functors. First we discuss in \seref{actions} actions of a set of natural transformations on a category. In \seref{bicat} we show how an adjoint pair in any bicategory can be formulated as a pair of invertible elements for a certain Morita context. Combining the notion of a Morita context over a ring with local units with the action of a class of natural transformations on a category, we then introduce the notion of locally adjoint functors in \seref{locadjfun}.
Following the philosophy of \cite{CMZ}, Frobenius properties of corings are related to the adjunction properties of the induction functor of a coring. For this reason, we study
in \seref{inductionadjoint} the induction functor $-\ot_A\cc:\Mm_A\to\Mm^\cc$ for an $A$-coring $\cc$, and we describe in \seref{nattranf} all the natural transformations from this functor to its left and right adjoint. In \seref{yoneda} we describe natural transformations between representable functors that are involved in the description of the Frobenius property as in \cite{CR, mio}. In order to allow a description of the quasi-Frobenius property, we repeat these procedures in a more general setting in \seref{coproduct}, involving a coproduct functor.
In \seref{frob} we introduce the notion of a locally Frobenius coring and a locally quasi-Frobenius coring, that coincides with the notion of a co-Frobenius coring, respectively quasi-co-Frobenius coring, if the base ring is a $PF$-ring. We give a characterization of locally quasi-Frobenius corings and prove some properties: we show that they provide examples of quasi-co-Frobenius corings over arbitrary base rings, they are locally projective as a left and right module over the base algebra (\thref{denseB}), and they are semiperfect if the base algebra is a QF-ring (\prref{semiperfect}).

Finally we apply all the obtained results in Sections \ref{se:frobenius} and \ref{se:QcF}, where we give a characterization of co-Frobenius and quasi-co-Frobenius corings and recover old characterizations of Frobenius corings.

\section{Preliminaries}\selabel{preliminaries}
Throughout this paper, $k$ will be a commutative base rings. All rings that we will consider
will be $k$-algebras, and categories will usually be $k$-linear. Unless otherwise stated,
functors will be covariant. For an object $X$ in a category $\Cc$, $X$ will also be our notation
for the identity morphism on $X$.
Let $I$ be any index set and $M$ an object in a category with products and coproducts.
We will denote $M^{(I)}$ for the coproduct (direct sum) and $M^I$ for the product. 

Let $R$ be a ring, possibly without unit. $\widetilde{\Mm}_R$ will denote the category of
right $R$-modules. For a ring $R$ with unit, $\Mm_R$ will denote the category of
unital right $R$-modules.

\subsection{Adjoint functors}
Let $\Cc$ and $\Dd$ be two categories and $F:\ \Cc \to \Dd$ and $G:\ \Dd\to \Cc$ two functors. We call $F$ a left adjoint of $G$, $G$ a right adjoint of $F$ or $(F,G)$ a pair of adjoint functors if and only if there exist natural isomorphisms 
\begin{equation}\eqlabel{thetaadjoint}
\theta_{C,D} :\ \Hom_\Dd(FC,D)\cong\Hom_\Cc(C,GD),
\end{equation}
for all $C\in\Cc$ and $D\in\Dd$.
This is equivalent to the existence of natural transformations
$\eta \in\Nat(\id_\Cc,GF)$ and $\varepsilon \in\Nat(FG,\id_\Dd)$, such that
\begin{eqnarray}
\varepsilon_{FC}\circ F(\eta_C)&=&FC,\qquad \forall C\in\Cc; \eqlabel{triang1}\\
G(\varepsilon_D)\circ\eta_{GD}&=&GD, \qquad \forall D\in\Dd. \eqlabel{triang2}
\end{eqnarray}

\subsection{Rings and corings}\selabel{ringscorings}

Let $A$ be a $k$-algebra. An \emph{$A$-ring} $(R,\mu,\eta)$ is an algebra (or monoid) in the monoidal category ${_A\Mm_A}$ consisting of $A\hbox{-}A$ bimodules and $A\hbox{-}A$ bilinear maps. There exists a bijective correspondence between a $A$-rings $R$ and 
$k$-algebras $R$ together with an algebra morphism $\eta:\ A\to R$.

The dual notion of an $A$-ring is an \emph{$A$-coring}, i.e.\ an $A$-coring $(\cc,\Delta_\cc,\varepsilon_\cc)$ is a coalgebra (or comonoid) in ${_A\Mm_A}$. Explicitly, an $A$-coring consists of an $A\hbox{-}A$ bimodule $\cc$ and two $A\hbox{-}A$ bilinear maps $\Delta_\cc:\cc\to \cc\ot_A\cc$
(the comultiplication) and $\varepsilon_\cc:\cc\to A$ (the counit), such that $(\Delta_\cc\otimes_A\cc)\circ \Delta_\cc=
(\cc\ot_A\Delta_\cc)\circ \Delta_\cc$ and 
$(\varepsilon_\cc\otimes_A\cc)\circ \Delta_\cc=(\cc\ot_A\varepsilon_\cc)\circ \Delta_\cc=
\Delta_\Cc$. For the comultiplication, we use the Sweedler-Heyneman notation, 
namely $\Delta(c)=c_{(1)}\otimes_Ac_{(2)}$ (summation understood) and $(\cc\otimes_A\Delta_\cc)\circ\Delta(c)
=(\Delta_\cc\otimes_A\cc)\circ\Delta(c)=c_{(1)}\otimes_Ac_{(2)}\otimes_Ac_{(3)}$. 

The category of right (resp. left) comodules over $\cc$ will be denoted by $\Mm^\cc$ (resp. ${^\cc\Mm}$). Recall that a right $\cc$-comodule $(M,\rho_M)$ consists of a right $A$-module $M$ and a right $A$-module map $\rho_M:M\to M\ot_A\cc$, $\rho_M(m)=m_{[0]}\ot_Am_{[1]}$, which satisfies the usual coassociativity and counit conditions.

For more details about the general theory of corings and their comodules, we refer the monograph \cite{BrzWis:book}.

The following elementary results from module theory will turn out to be useful in the sequel.

\begin{lemma}\lelabel{sumproduct}
Let $A$ and $B$ be objects in a category $\Aa$, and $I$ an index set.
If $A^{(I)}$ and $B^I$ exist, then
$\Hom(A^{(I)},B)\cong\Hom(A,B^I)\cong (\Hom(A,B))^I$.
\end{lemma}

\begin{proof}
For $\ell\in I$, let $\iota_\ell:\ A\to A^{(I)}$  
and $\pi_\ell:B^I\to B$ be the canonical 
canonical coproduct and product maps. Consider the diagram
\[\xymatrix{
A \ar[rr]^-{f^\circ} \ar[d]_{\iota_\ell} \ar[drr]^{f_\ell} && B^I \ar[d]^{\pi_\ell}  \\
A^{(I)} \ar[rr]_{f_\circ} && B
}\]
Any morphism $f^\circ \in\Hom(A,B^I)$ as well as any morphism $f_\circ\in\Hom(A^{(I)},B)$ is completely determined by the family of morphisms $f_\ell: A\to B$, $\ell\in I$. 
\end{proof}

\begin{lemma}\lelabel{M^B}
Consider a ring morphism $B\to A$ and take any $M\in{_A\Mm_B}$, then ${_A\Hom_B}(A,M)\stackrel{\xi_M}{\cong} M^B:=\{m\in M~|~ bm=mb, \textrm{ for all } b\in B\}$.
\end{lemma}

\begin{proof}
For any $f\in{_A\Hom_B}(A,M)$, we obtain $bf(1_A)=f(b)=f(1_A)b$, consequently $f(1_A)\in M^B$. Conversely for any $x\in M^B$, define $f_x\in{_A\Hom_B}(A,M)$ as $f_x(a)=ax$. One can easily check that this correspondence is bijective.
\end{proof}

\subsection{Local units and local projectivity}\selabel{locproj}
Let $R$ be a non-unital $B$-ring and $(M,\mu_M)$ a right $R$-module, i.e.\ $M$ is a right $B$-module and $\mu_M:M\ot_BR\to M$ is an associative right $B$-linear multiplication map. We say that $R$ has right local units on $M$ if for every finitely generated right $B$-submodule $N$ of $M$, there exists an element $e\in R^B$ such that $n\cdot e=n$ for all $n\in N$. We call $e$ a (right) local unit for $N$. One can easily prove that $R$ has right local units on $M$ if and only if $R$ has right local units on every singleton $\{m\}\subset M$.
We say that $R$ is a ring with right local units if $R$ has right local units on $R$, where we consider the regular right $R$-module structure on $R$. 
The following Theorem can be viewed as a structure theorem for modules over rings with local units, and should be compared to similar results for rings with idempotent local units (see \cite[Lemma 2.10]{V:locunit} and \cite[Lemma 4.3]{CDV:colimit}).

\begin{theorem}\thlabel{locunit}
Let $R$ be a $B$-ring (without unit). Let $\Mm$ be a full subcategory of $\widetilde{\Mm}_R$.
Then $R$ has local units on all $M\in \Mm$ if and only if there exists a full subcategory $\Nn$
of $\widetilde{\Mm}_B$, such that every $M\in \Mm$ is generated by objects of $\Nn$ as a right $B$-module,
 and, for all $N\in \Nn$ and $f\in\Hom(N,M)$, we can find an $e\in R^B\cong{_B\Hom_B}(B,R)$ such that $f=f_e\circ f$, where 
\[
\xymatrix{
f_e=\mu_M\circ (M\ot_B e): M \cong M\ot_BB \ar[rr]^-{M\ot_Be} && M\ot_BR \ar[rr]^-{\mu_M} && M.
}
\]
In other words, $\Mm$ is generated by $B$-modules on which there exists a local unit.
\end{theorem}

\begin{proof}
Suppose first that the subcategory $\Nn$ exists. Take any $M\in\Mm$. Since $\Nn$ generates $\Mm$, we can find a family of right $B$-modules $(N_i)_{i\in I}$ in $\Nn$ such that there exists a surjective map $\pi : \coprod_{i\in I} N_i \to M$. Consequently, for any $m\in M$, we can write $m=\sum_{i\in J}\pi(n_i)$ where $J$ is a finite subset of $I$. We show by induction on the cardinality of $J$ that we can find a local unit $e\in R^B$. If the cardinality of $J$ equals one, then $m=\pi(n)$ for some $n\in N$. 
We know that there exists an element $e\in R^B$ such that $\pi=f_e\circ\pi$.
Consequently $\pi(n)e=f_e\circ\pi(n)=\pi(n)$, so $e$ is a unit for $\pi(n)=m$. Now suppose $m=\sum_{i=1}^k\pi(n_i)$ with $n_i\in N_i$ and $k> 1$. By the induction hypothesis we can find a local unit $e\in R^B$ for $\sum_{i=1}^{k-1}\pi(n_i)$ and a local unit $e'\in R^B$ for $\pi(n_k)-\pi(n_k)e$. Then $e''=e+e'-ee'$ 
is a local unit for $m$ since
\begin{eqnarray*}
m e'' &=& \left(\sum_{i=1}^k\pi(n_i)\right)(e+e'-ee')
=\left(\sum_{i=1}^{k-1}\pi(n_i)\right)(e+e'-ee')+\pi(n_k)(e+e'-ee')\\
&=&\sum_{i=1}^{k-1}\pi(n_i)+\sum_{i=1}^{k-1}\pi(n_i)e'-\sum_{i=1}^{k-1}\pi(n_i)e'
+\pi(n_k)e+(\pi(n_k)-\pi(n_k)e)e'\\
&=&\sum_{i=1}^{k-1}\pi(n_i)+\pi(n_k)=m.
\end{eqnarray*}

Conversely, let $\Mm$ be a subcategory of $\widetilde{\Mm}_R$ on which $R$ has local units. We define $\Nn$ as the category consisting of finitely generated $B$-submodules of modules in $\Mm$. Then clearly $\Nn$ generates $\Mm$ in $\Mm_B$.
Moreover, 
for any $N\in\Nn$, $M\in\Mm$ and $f\in \Hom_B(M,N)$,
we obtain that $\im f$ is a finitely generated $B$-submodule of $M$. By the definition of a module with local units, we can find a local unit $e\in R^B$ for $\im f$. Consequently, $f_e\circ f(n)=f(n)e=f(n)$ for all $n\in N$.
\end{proof}

Recall from \cite{ZH} that a right $A$-module is called locally projective if for any commutative diagram
in $\Mm_A$ with exact rows
\[
\xymatrix{
0 \ar[r] & F \ar[r]^i & M \ar[d]^g \\
& N' \ar[r]_f & N \ar[r] & 0
}
\]
where $F$ is finitely generated, there exists a right $A$-linear map 
$h:\ M\to N'$ such that $g\circ i = f\circ h\circ i$.
In \cite{Gar} it is shown that $M$ is locally projective if and only if for any finitely generated $A$-submodule $F\subset M$, there exists a finite dual basis $\{e_i,f_i\}\subset M\times M^*$. More generally, a $B\hbox{-}A$ bimodule $M$ is called $\Rr$-locally projective for an additive subset $\Rr\subset M^*=\Hom_A(M,A)$ if for any finite subset $N$ of $M$ we can find a finite set $\{e_i,f_i\}\subset M\times\Rr$ such that $e_if_i(n)=n$ for all $n\in N$ and $e_if_i(bm)=be_if_i(m)$ for all $m\in M$. Obviously, $\Rr$-local projectivity implies $\Ss$-local projectivity if $\Rr\subset\Ss$. In particular $\Rr$-local projectivity implies local projectivity.
The relationship between local projectivity and local units is discussed in general in \cite{V:locunit}. 

\subsection{Modules versus comodules}\selabel{modversuscomod}
Let $\cc$ be an $A$-coring. It is well-known that there exists an adjunction $(\Ff^\cc,\Gg^\cc)$ between the forgetful functor $\Ff^\cc:\Mm^\cc\to \Mm_A$ and the induction functor  $\Gg^\cc=-\ot_A\cc:\Mm_A\to\Mm^\cc$ (see \seref{induction}); this implies that we have
a natural isomorphism $\Hom^\cc(M,N\ot_A\cc)\cong\Hom_A(M,N)$ for all $M\in\Mm^\cc$ and $N\in\Mm_A$. Consequently, $\cc^*=\Hom_A(\cc,A)\cong\End^\cc(\cc)$ is an $A$-ring with unit $\varepsilon_\cc$ and multiplication given by
$$f*g(c)=f(g(c_{(1)})c_{(2)}),$$
for all $f, g\in \cc^*$ and $c\in \cc$. In a similar way, the left dual $\*c={_A\Hom}(\cc,A)$ is an $A$-ring with unit $\varepsilon_\cc$ and multiplication
$$f*g(c)=g(c_{(1)}f(c_{(2)})).$$
Finally, $\*c^*={_A\Hom_A}(\cc,A)$ is a $k$-algebra with unit $\varepsilon_\cc$ and multiplication
$$f*g(c)=g(c_{(1)})f(c_{(2)}).$$
Note that if $A$ is commutative and $\cc$ is an $A$-coalgebra (i.e.\ the left and right $A$-action on $\cc$ coincide), then $\*c=\cc^*=\*c^*$ with the opposite multiplication (see e.g.\ \cite{DNR}).
Furthermore, every right comodule has a right  $\*c$-module structure, given by 
$$m\cdot f=m_{[0]}f(m_{[1]}),$$
for any $m\in M\in \Mm^\cc$ and $f\in\*c$. In a similar way, every left $\cc$-comodule has a left $\cc^*$-module structure.

Let $\Rr$ be any additive subset of $\*c$ such that $A\Rr\subset \Rr A$ and take any $M\in\Mm_{\*c}$. The $\Rr$-rational part of $M$ is defined as
$$\Rat_\Rr(M)=\{m\in M~|~\exists~m_i\in M, c_i\in\cc, \textrm{ such that } \sum_i m\cdot r=m_ir(c_i), \forall r\in\Rr \}.$$
We say that $M$ is $\Rr$-rational if $\Rat_\Rr(M)=M$. Denote by $T$ the subring of $\*c$, generated by $A$ and $\Rr$. It was proved in \cite{CVW2} (see also \cite{Abu}) that if $\cc$ is $\Rr$-locally projective as a left $A$-module, then every $\Rr$-rational $\Rr$-faithful $\*c$-module is a $\cc$-comodule with coaction given by
$\rho(m)=\sum_im_i\ot_Ac_i$ if and only if $m\cdot r=\sum_im_ir(c_i)$, for all $r\in\Rr$.
This defines a functor 
$$\Rat_\Rr : \Mm_{\*c}\to \Mm^\cc.$$
The category of all $\Rr$-rational $\Rr$-faithful right $\*c$-comodules and $T$-linear maps is isomorphic to the category of right $\cc$-comodules.
The coring $\cc$ is finitely generated and projective as a left $A$-module, if and only if $\Rat_\Rr(\*c)=\*c$.

\subsection{Frobenius corings}

An $A$-coring is called \emph{Frobenius} if $\cc$ and $\*c$ are isomorphic as $A\hbox{-}\*c$ bimodules. It is well-known that this notion is left-right symmetric: $\cc$ is Frobenius if and only if $\cc$ and $\cc^*$ are isomorphic as $\cc^*$-$A$ bimodules (see also \coref{Frobenius}). 

The notion of co-Frobenius and quasi-co-Frobenius coalgebra can in a natural way 
be generalized to the setting of corings, as follows:

\begin{Definition}
An $A$-coring is called \emph{left co-Frobenius} if and only if there exists an
$A\hbox{-}\*c$ bimodule monomorphism $j:\ \cc\to \*c$. 

$\cc$ is called \emph{left quasi-co-Frobenius} if there exists an $A\hbox{-}\*c$ bimodule monomorphism 
$j:\ \cc\to(\*c)^{I}$ for some index set $I$.

Right co-Frobenius and quasi-co-Frobenius corings can be introduced in a similar way, replacing $\*c$
by $\cc^*$ and requiring the existence of $\cc^*$-$A$ bimodule monomorphisms.
\end{Definition}

\begin{remark}
A left quasi-co-Frobenius coalgebra over a field $k$ is usually defined as a $k$-coalgebra $C$ such that there exists a monomorphism  $j:\ C\to (C^*)^{(I)}$ of left $C^*$-modules. Considering $C$ as a $k$-coring,
we remark first that the convolution product on $C^*$  is opposite to the multiplication in $\*c$,
if we use convention introduced in \seref{modversuscomod}. Secondly, it was proved in \cite[Theorem 1.3]{GTN} that the existence of a left $C^*$-linear monomorphism $j:\ C\to (C^*)^{(I)}$ is equivalent to the existence of a left $C^*$-linear monomorphism $j':\ C\to (C^*)^I$.
\end{remark}

Let $(F,G)$ be a pair of adjoint functors. The pair $(F,G)$ is called a \emph{Frobenius pair} of functors, if $G$ is also a left adjoint of $F$.

\subsection{Morita contexts}\selabel{intrMorita}

Recall (see \cite[Remarks p 389, Examples 1.2]{Mul}, \cite[Remark 3.2]{BohmBrz:cleft}) that a Morita context can be identified with a $k$-linear category with two objects $a$ and $b$. The algebras of the Morita contexts are $\End(a)$ and $\End(b)$, the connecting bimodules are $\Hom(a,b)$ and $\Hom(b,a)$ and multiplication and bimodule maps are given by composition. We denote this context as follows
$$\NN(a,b)=(\End(a),\End(b),\Hom(b,a),\Hom(a,b),\circ,\bul).$$
This can be summarized by the following diagram.
\[
\xymatrix{
a \ar@(dl,ul)^{\End(a)} \ar@/^/[rr]^{\Hom(a,b)} && b \ar@/^/[ll]^{\Hom(b,a)} \ar@(dr,ur)_{\End(b)} 
}
\]
If $j\in\Hom(a,b)$ and $\barj\in\Hom(b,a)$ are such that $\barj\circ j = a$ and $j\circ\barj=b$, then we
call $(j,\barj)$ a pair of \emph{invertible elements}. This means that $j$ and $\barj$ are inverse
isomorphisms between $a$ and $b$. As a Morita context $\MM$ with an invertible pair $(j,\barj)$ is always strict, we say that $\MM$ is strict \emph{by} $(j,\barj)$.

A morphism of Morita contexts $$\mmm: \MM=(A,B,P,Q,\mu,\tau) \to \MM'=(A',B',P',Q',\mu',\tau')$$ consists of two algebra maps $\mmm_1:A\to A'$ and $\mmm_2:B\to B'$, an $A'$-$B'$ bimodule map $\mmm_3:P\to P'$ and a $B'$-$A'$ bimodule map $\mmm_4:Q\to Q'$ such that $\mmm_1\circ\mu=\mu'\circ(\mmm_3\ot_B\mmm_4)$ and $\mmm_2\circ\tau=\tau'(\mmm_4\ot_A\mmm_3)$.

There are two canonical ways to construct new Morita contexts out of an existing one, without adding or removing any information.
\begin{enumerate}[(i)]
\item
The \emph{opposite} of a Morita context $\MM=(A,B,P,Q,\mu,\tau)$ is the Morita context $\MM^{\rm op}=(A^{\rm op},B^{\rm op},Q,P,\mu^{\rm op},\tau^{\rm op})$, where $\mu^{\rm op}(q\ot_{B^{op}} p)=\mu(p\ot_B q)$ and $\tau^{\rm op}(p\ot_{A^{op}} q)=\tau(q\ot_A p)$. 
An anti-morphism of Morita contexts $\mmm : \MM\to \MM'$ is a morphism $\mmm:\MM\to {\MM'}^{\rm op}$. It consists of two algebra maps $\mmm_1:\ A\to {A'}^{\rm op}$ and $\mmm_2:\ B\to {B'}^{\rm op}$, an $A'$-$B'$ bimodule map $\mmm_3:Q\to P'$ and a $B'$-$A'$ bimodule map $\mmm_4:P\to Q'$ such that $\mmm_1\circ\mu={\mu'}^{\rm op}\circ(\mmm_4\ot_B\mmm_3)$ and $\mmm_2\circ\tau={\tau'}^{\rm op}(\mmm_3\ot_A\mmm_4)$.
\item
The \emph{twisted} of a Morita context $\MM=(A,B,P,Q,\mu,\tau)$ is the Morita context $\MM^{\rm t}=(B,A,Q,P,\tau,\mu)$. 
\end{enumerate}

\section{Locally adjoint functors}\selabel{adjoint}\selabel{locadj}

\subsection{Action of a set of natural transformations on a category}\selabel{actions}

Let $F:\Cc\to\Dd$ be a functor and consider 
a semigroup of natural transformations $\Phi\subset\Nat(F,F)$. 
We define for all $\alpha\in\Phi$, $C,C'\in \Cc$ and $f:F(C)\to F(C')$ in $\Dd$,
$$\alpha\cdot f:=\alpha_{C'}\circ f:F(C)\to F(C').$$
This defines an action of $\Phi$ on $\Hom_\Dd(F(C),F(C'))$, indeed
for $\alpha,\beta\in\Phi$, $C,C'\in \Cc$ and $f:F(C)\to F(C')$ in $\Dd$, we have
$$ (\alpha\circ\beta)\cdot f =
(\alpha\circ\beta)_{C'} \circ f = \alpha_{C'}\circ \beta_{C'} \circ f = \alpha\cdot (\beta\cdot f).$$
Since this action exists for all choices of $C,C'\in \Cc$, we will say that $\Phi$ acts on $\Cc$.
We say that $\Phi$ acts unital on $\Cc$, if there exists an element $e\in \Phi$ such that for all $f:F(C)\to F(C')$ in $\Dd$ with $C,C'\in \Cc$, we have $e\cdot f=f$. We say that $\Phi$ acts with local units on $\Cc$ if an only if there exists a generating subcategory $\Ee\subset \Cc$, such that for all $f:F(E)\to F(C)$ in $\Dd$ with $C\in \Cc$ and $E\in \Ee$, there exists an $e^E\in\Phi$ such that $e^E\cdot f=f$.

\begin{example}[\textbf{action of a ring with local units}]
Let $R$ be a (non-unital) $B$-ring, and $\Mm$ be a full subcategory of $\widetilde{\Mm}_R$. Consider the forgetful functor $U:\Mm\to\Mm_B$. Then we have a map $R^B\to\Nat(U,U)$. Indeed take any $M\in\Mm$, $e\in R^B$ and define $\alpha^e_M(m)=m\cdot e$ for all $m\in M$, then $\alpha^e_M$ is a right $B$-linear map that is natural in $M$. Now consider a subcategory $\Nn\subset\Mm_B$. Following the procedure of this section, $R^B$, viewed as a subset of $\Nat(U,U)$, has the following action on $\Nn$. For any right $B$-linear map $f:N\to M$ with $N\in \Nn$ and $M\in \Mm$ we define $e\cdot f=\alpha^e_M\circ f$. Suppose now that $\Nn$ generates $\Mm$ in $\Mm_B$. Then it follows from \thref{locunit} that $R$ acts with local units on all objects in $\Mm$ if and only if $R^B$ (as natural transformations) acts unitally on $\Nn$, in other words, if and only if $R^B$ (as natural transformations) acts with local units on $\Mm$.
\end{example}

Consider now a functor $G:\Dd\to \Cc$ and let $\Gamma$ be a semigroup of natural transformations $\Gamma\subset\Nat(G,G)$. We define for all $\alpha\in\Gamma$ and $f:C\to G(D)$ with $C\in \Cc$ and $D\in\Dd$,
$$\alpha\cdot f:=\alpha_D\circ f:C\to G(D).$$
One can easily check that defines an action of $\Gamma$ on $\Hom_\Cc(C,G(D))$. Since this action exists for all choices of $C$, we will say that $\Gamma$ acts on $\Cc$.
We say that $\Gamma$ acts unital on $\Cc$, if there exists an element $e\in \Gamma$ such that for all $f:C\to G(D)$ in $\Cc$ with $C\in \Cc$ and $D\in\Dd$, we have $e\cdot f=f$. We say that $\Gamma$ acts with local units on $\Cc$ if an only if there exists a generating subcategory $\Ee\subset \Cc$, such that for all $f:E\to G(D)$ in $\Cc$ with $E\in \Ee$ and $D\in\Dd$, there exists an $e^E\in\Gamma$ such that $e^E\cdot f=f$.

\begin{example}
Let $R$ be a ring (with unit), and $\Mm_R$ the category of right $R$-modules. Denote by $\bf 1$ the category with a unique object $*$ and a unique (endo)morphism. Consider a functor $G:{\bf 1}\to \Mm_R$. Then $G$ is completely determined by $G(*)=M$, and furthermore $\Nat(G,G)\cong\End_R(M)$. Therefore, for any semigroup $\Gamma$ of endomorphisms of $M$, we have a natural map $\Gamma\to\Nat(G,G)$. We know that $R$ is a generator for $\Mm_R$ and $\Hom_R(R,G(*))=\Hom_R(R,M)\cong M$. The action $\Gamma$ (considered as set of natural transformations) on $M\in\Mm_R$ coincides with the canonical action of $\Gamma$ (considered as set of endomorphisms) on $M$.
\end{example}

\subsection{Adjoint functors and Morita contexts}\selabel{bicat}

The notion of a bicategory was introduced in \cite{Ben}, see also e.g.\ \cite[Chapter XII]{McLane}. We will use the notion of a $\Vv$-enriched bicategory, where $\Vv$ is a monoidal category.
A $\Vv$-enriched bicategory $\Bb$ consists of the following data,
\begin{enumerate}[(i)]
\item a class of objects $A,B,\ldots$ which are called $0$-cells;
\item for every two objects $A$ and $B$, a $\Vv$-enriched category $\Hom(A,B)$,
whose class of objects, which are called $1$-cells, we denote by $\Hom_1(A,B)$. 
We write $f:A\to B$ for a $1$-cell $f\in\Hom_1(A,B)$. 
The set of morphisms between two
$1$-cells $f, g\in\Hom_1(A,B)$ is denoted by ${^A\Hom_2^B}(f,g)\in\Vv$. We call these morphisms $2$-cells and denote them as $\alpha:f\rightarrow g$;
\item[]\hspace{-.9cm}furthermore there exist compositions $\bul, \ubul$ and $\circ$ as follows
\item for all $f\in\Hom_1(A,B)$ and $g\in\Hom_1(B,C)$, we have $f\bul_B g\in\Hom_1(A,C)$;
\item for all 
$\alpha\in{^A\Hom_2^B}(f,g)$ and $\beta\in{^B\Hom_2^C}(h,k)$, we have
$\alpha \ubul_B\beta:f\bul_B h\rightarrow g\bul_B k$; by the $\Vv$-enriched property this can be expressed as
$$\ubul:{^A\Hom_2^B}(f,g)\ot{^B\Hom_2^C}(h,k)\to {^A\Hom_2^C}(f\bul_Bh,g\bul_B k);$$
\item for all $f,g,h\in\Hom_1(A,B)$ such that $\alpha:f\rightarrow g$ and $\beta:g\rightarrow h$, then $\beta\circ \alpha:f\rightarrow h$, this is just the composition of morphisms in the category $\Hom(A,B)$.
\item For all $0$-cells $A$ in $\Bb$ there exists a $1$-cell $\id_A:A\to A$ such that $\id_A\bul_A f\cong f$ and $g\bul_A\id_A\cong g$ for all $0$-cells $B$ and all $f\in\Hom_1(A,B)$ and $g\in\Hom_1(B,A)$.
\end{enumerate}
For all compatibility conditions we refer to \cite{Ben}. Let us recall the interchange law
\begin{equation}\eqlabel{interchange}
(\alpha\ubul\beta)\circ(\gamma\ubul\delta)=(\alpha\circ\gamma)\ubul(\beta\circ\delta),
\end{equation}
for $\alpha\in{^A\Hom_2^B}(a,c)$, $\beta\in{^B\Hom_2^C}(b,d)$, $\gamma\in{^A\Hom_2^B}(c,e)$ and $\delta\in{^B\Hom_2^C}(d,f)$.
From \equref{interchange} we immediately deduce that
\begin{equation}\eqlabel{interchange2}
(\alpha\ubul b)\circ(c\ubul \beta)=\alpha\ubul\beta=(a\ubul\beta)\circ(\alpha\ubul d),
\end{equation}
for all $\alpha\in{^A\Hom_2^B}(a,c)$ and $\beta\in{^B\Hom_2^C}(b,d)$.

Recall that by the Coherence Theorem (see \cite{McLanePar} and \cite{Str:book}), we are allowed to preform calculations for bicategories in the simpler formalism of $2$-categories, i.e.\ such that $\id_A=A$ and the isomorphisms
$m\bul_AA\cong m$, $(m\bul_A n)\bul_B p\cong m\bul_A(n\bul_B p)$ and the corresponding isomorphisms for $\ubul$ are the identity morphisms. The basic example of such a $2$-category is $\CAT$, the bicategory consisting of categories, functors and natural transformations (see \seref{locadjfun}).

Morita theory can be developed naturally within the framework of bicategories. A \emph{Morita context} in a bicategory $\Bb$ is a sextuple $(A,B,p,q,\mu,\tau)$, where $A$ and $B$ are $0$-cells, $p\in\Hom_1(A,B)$, $q\in\Hom_1(B,A)$, $\mu\in{^A\Hom_2^A}(p\bul_Bq,A)$ and $\tau\in{^B\Hom_2^B}(q\bul_Ap,B)$ such that $q\ubul_A \mu=\tau\ubul_B q$ and $p\ubul_B\tau=\mu\ubul_Ap$.

An \emph{adjoint pair} in $\Bb$ is a sextuple $(A,B,p,q,\mu,\nu)$, where $A$ and $B$ are $0$-cells, $p\in\Hom_1(A,B)$, $q\in\Hom_1(B,A)$, $\mu\in{^A\Hom_2^A}(p\bul_Bq,A)$ and $\nu\in{^B\Hom_2^B}(B,q\bul_Ap)$, such that $(\mu\ubul_Ap)\circ(p\ubul_B\nu)=p$ and $(q\ubul_A\mu)\circ(\nu\ubul_Aq)=q$.

Let $\Vv$ be a monoidal category with coequalizers, this is a monoidal category that possesses coequalizers and in which the tensor product preserves these coequalizers, see \cite{Bohm:internal}.
Then we can construct a $\Vv$-enriched bicategory $\Bim(\Vv)$ as follows. 
\begin{itemize}
\item $0$-cells are the algebras in $\Vv$;
\item $1$-cells are bimodules between those algebras;
\item $2$-cells are bimodule maps;
\item the composition of an $A$-$B$ bimodule $M$ and a $B$-$C$ bimodule $N$ is given by the following coequalizer $N\ot_B N$
\[
\xymatrix{
M\ot B\ot N \ar@<.5ex>[r] \ar@<-.5ex>[r] & M\ot N \ar[r] & M\ot_B N.
}\]
\end{itemize}
In the situation where $\Vv$ is the category of abelian groups, $\Bim(\Vv)=\Bim$ is the bicategory
of rings, bimodules and bimodule maps.

\begin{theorem}\thlabel{Morita}
Let $A$ and $B$ be two $0$-cells in a $\Vv$-enriched bicategory $\Bb$,
and $p: A\to B$, $q:B\to A$ two $1$-cells. Consider $Q={^B\Hom_2^A}(q,q)$,
$P={^A\Hom_2^B}(p,p)^{\rm op}$, $N={^B\Hom_2^B}(B,q\bul_A p)$  
and $M={^A\Hom_2^A}(p\bul_B q,A)$.
\begin{enumerate}
\item 
$Q$ and $P$ are $0$-cells in $\Bim(\Vv)$;
\item 
$N$ is a $1$-cell from $Q$ to $P$;
\item 
$M$ is a $1$-cell from $P$ to $Q$;
\item there exist two maps $\newtriangle:N\otimes_PM\to Q$ and $\newblacktriangle:M\otimes_QN\to P$;
\item $\MM(p,q)=(Q,P,N,M,\newtriangle,\newblacktriangle)$ is a Morita context in $\Bim(\Vv)$.
\end{enumerate}
\end{theorem}

\begin{proof}
\ul{(1)}. The `vertical' composition $\circ$ of $2$-cells in $\Bb$ defines an associative multiplication on $Q={^B\Hom_2^A}(q,q)$ and  ${^A\Hom_2^B}(p,p)$
with units $q\in Q={^B\Hom_2^A}(q,q)$ and $p\in{^A\Hom_2^B}(p,p)$. We will denote $\ast$ for the opposite multiplication in ${^A\Hom_2^B}(p,p)$.\\
\ul{(2)}. Take $\alpha,\alpha'\in Q$, $\beta,\beta'\in P$ and $\gamma\in N$. We define $\alpha\cdot \gamma=(\alpha\ubul_A p)\circ\gamma$ and $\gamma\cdot\beta=(q\ubul_A \beta)\circ\gamma$.
Then we find
\begin{eqnarray*}
&&\hspace*{-2cm} 
(\alpha \cdot \gamma)\cdot \beta=((\alpha\ubul_A p)\circ\gamma)\cdot \beta 
=(q\ubul_A \beta)\circ((\alpha\ubul_A p)\circ\gamma)\\
&\equal{\equref{interchange2}}&(\alpha\ubul_A p)\circ((q\ubul_A \beta)\circ\gamma)
=\alpha\cdot(\gamma\cdot \beta).
\end{eqnarray*}
Both actions are associative:
\begin{eqnarray*}
&&\hspace*{-2cm} 
\alpha'\cdot(\alpha\cdot \gamma)=\alpha'\cdot((\alpha\ubul_A p)\circ\gamma)
=(\alpha'\ubul_A p)\circ((\alpha\ubul_A p)\circ\gamma)\\
&\equal{\equref{interchange}}&((\alpha'\circ\alpha)\ubul_A p)\circ\gamma)
=(\alpha'\circ\alpha)\cdot\gamma;\\
&&\hspace*{-2cm} 
(\gamma\cdot\beta)\cdot\beta'=((q\ubul_A \beta)\circ\gamma)\cdot\beta'
=(q\ubul_A \beta')\circ((q\ubul_A \beta)\circ\gamma)\\
&\equal{\equref{interchange}}&(q\ubul_A (\beta'\circ \beta))\circ\gamma
=\gamma\cdot(\beta'\circ\beta)=\gamma\cdot(\beta\ast\beta').
\end{eqnarray*}
Obviously, $q\in Q$ and $p\in P$ act trivially on $N$.\\
\ul{(3)}. This statement is dual to $(2)$. We only give the definition of the actions and leave further verification to the reader. Take $\alpha\in Q$, $\beta\in P$ and $\delta\in M$, then
$\beta\cdot\delta=\delta\circ(\beta\ubul_Bq)$ and $\delta\cdot\alpha=\delta\circ(p\ubul_B\alpha)$.\\
\ul{(4)}. Take $\gamma\in N$ and $\delta\in M$. Then we define $\gamma\newtriangled\delta=(q\ubul_A\delta)\circ(\gamma\ubul_Bq)$ and $\delta\newblacktriangled\gamma=(\delta\ubul_Ap)\circ(p\ubul_B\gamma)$:
\begin{eqnarray*}
\xymatrix{
\gamma\newtriangled\delta :~~
q\cong B\bul_Bq \ar[r]^-{\gamma\ubul_Bq} & q\bul_Ap\bul_Bq \ar[r]^-{q\ubul_A\delta} & q\bul_AA \cong q;}
\\
\xymatrix{
\delta\newblacktriangled\gamma :~~
p\cong p\bul_BB \ar[r]^-{p\ubul_B\gamma} & p\bul_Bq\bul_Ap \ar[r]^-{\delta\ubul_Ap} & A\bul_Ap \cong p.
}
\end{eqnarray*}
Let us check that $\newtriangle$ is $P$-balanced and that $\newblacktriangle$ is $Q$-balanced.
For $\beta\in P$ and $\alpha\in Q$, we compute that
\begin{eqnarray*}
&&\hspace*{-2cm}
(\gamma\cdot \beta)\newtriangled\delta = (q\ubul_A\delta)\circ((\gamma\cdot\beta)\ubul_Bq)
= (q\ubul_A\delta)\circ(((q\ubul_A \beta)\circ\gamma)\ubul_Bq)\\
&=& (q\ubul_A\delta)\circ(q\ubul_A \beta\ubul_Bq)\circ(\gamma\ubul_Bq)
=(q\ubul_A(\delta\circ(\beta\ubul_Bq)))\circ(\gamma\ubul_Bq)\\
&=&(q\ubul_A(\beta\cdot \delta))\circ(\gamma\ubul_Bq)
=\gamma\newtriangled(\beta\cdot \delta);\\
&&\hspace*{-2cm}
\delta\newblacktriangled(\alpha\cdot\gamma)
= (\delta\ubul_Ap)\circ(p\ubul_B(\alpha\cdot\gamma))
= (\delta\ubul_Ap)\circ(p\ubul_B((\alpha\ubul_A p)\circ\gamma))\\
&=& (\delta\ubul_Ap)\circ(p\ubul_B\alpha\ubul_A p)\circ(p\ubul_B\gamma)
= ((\delta\circ(p\ubul_B\alpha))\ubul_Ap)\circ(p\ubul_B\gamma)\\
&=& ((\delta\cdot\alpha)\ubul_Ap)\circ(p\ubul_B\gamma)
= (\delta\cdot \alpha)\newblacktriangled \gamma.
\end{eqnarray*}
\ul{(5)}. For $\gamma, \gamma'\in N$ and $\delta, \delta'\in M$, we compute
\begin{eqnarray*}
&&\hspace*{-2cm}
\gamma\cdot(\delta\newblacktriangled\gamma')
= (q\ubul_A(\delta\newblacktriangled\gamma'))\circ\gamma
= (q\ubul_A((\delta\ubul_Ap)\circ(p\ubul_B\gamma')))\circ\gamma\\
&=& (q\ubul_A\delta\ubul_Ap)\circ(q\ubul_Ap\ubul_B\gamma')\circ\gamma
= (q\ubul_A\delta\ubul_Ap)\circ\gamma\ubul_B\gamma'\\
&=& (q\ubul_A\delta\ubul_Ap)\circ(\gamma\ubul_Bq\ubul_Ap)\circ\gamma'
= (((q\ubul_A\delta)\circ(\gamma\ubul_Bq))\ubul_Ap)\circ\gamma'\\
&=& ((\gamma\newtriangled\delta)\ubul_Ap)\circ\gamma'
= (\gamma\newtriangled\delta)\cdot\gamma'.
\end{eqnarray*}
A similar computation shows that $\delta\cdot(\gamma\newtriangled\delta')=(\delta\newblacktriangled\gamma)\cdot\delta'$.
\end{proof}

Recall from \seref{intrMorita} that $(\mu,\nu)\in M\times N$ is a pair of invertible elements for $\MM(p,q)$
if $\newtriangle(\nu\ot_P\mu)=Q$ and $\newblacktriangle(\mu\ot_A \nu)=P$.
Obviously the Morita context $\MM(p,q)$ is strict if there exists a pair of invertible elements
(but not conversely).
Comparing the definitions of an adjoint pair and of $\newtriangle$ and $\newblacktriangle$,
we immediately obtain the following result.

\begin{theorem}\thlabel{Morita-adjoint}
Let $A$ and $B$ be two $0$-cells in $\Bb$, and $p: A\to B$ and $q:B\to A$ $1$-cells.
With notation as above, $(\mu,\nu)$ is a pair of invertible elements for $\MM(p,q)$
if and only $(A,B,p,q,\mu,\nu)$ is an adjoint pair in $\Bb$.
\end{theorem}

\subsection{Comatrix coring contexts}
Let $k$ be a commutative ring, and consider the bicategory
$\Bb=\Bim(\Mm_k)$.
In the literature (see \cite{BGT})
an adjoint pair in $\Bb$ is also termed a comatrix coring context. 
An interesting aspect of comatrix coring contexts is that they can be used to construct
certain corings, called comatrix corings.
A comatrix coring context is a sextuple $(A,B,\Sigma^\dagger,\Sigma,\varepsilon,\eta)$, where $A$ and $B$ are rings, $\Sigma^\dagger\in{_A\Mm_B}$, $\Sigma\in{_B\Mm_A}$, and $\varepsilon :\Sigma^\dagger\ot_B\Sigma\to A$ and $\eta:B\to \Sigma\ot_A\Sigma^\dagger$ 
are bimodule maps such that the following diagrams commute
\begin{equation}\eqlabel{ccc}
\xymatrix{
\Sigma \ar[rr]^-{\cong} \ar[d]_{\cong} && B\ot_B\Sigma \ar[d]^{\eta\ot_B\Sigma}\\
\Sigma\ot_AA && \Sigma\ot_A\Sigma^\dagger\ot_B\Sigma \ar[ll]^{\Sigma\ot_A\varepsilon}
}\qquad
\xymatrix{
\Sigma^\dagger \ar[rr]^-{\cong} \ar[d]_{\cong} && \Sigma^\dagger\ot_BB \ar[d]^{\Sigma^\dagger\ot_B\eta}\\
A\ot_A\Sigma^\dagger && \Sigma^\dagger\ot_B\Sigma\ot_A\Sigma^\dagger \ar[ll]^{\varepsilon\ot_A\Sigma^\dagger}
}
\end{equation}
The existence of a comatrix coring context implies that $\Sigma$ is finitely generated and projective as a right $A$-module and $\Sigma^\dagger\cong\Sigma^*$.

We can also consider comatrix coring contexts in the bicategory $\Frm(\Mm_k)$. The $0$-cells in
$\Frm(\Mm_k)$ are firm algebras, $1$-cells are firm bimodules, and $2$-cells are
bimodule maps. The existence of a comatrix coring context now implies that $\Sigma$ is
right $A$-firmly projective in the sense of \cite{V:equivalence}. Corings arising from 
comatrix coring contexts in $\Frm(\Mm_k)$ are known as infinite comatrix corings,
see \cite{CDV:colimit,V:equivalence}.

Consider $\Sigma^\dagger\in{_A\Mm_B}$ and $\Sigma\in{_B\Mm_A}$. Applying \thref{Morita},
we obtain a Morita context $\MM(\Sigma,\Sigma^\dagger)=(Q={_B\End_A}(\Sigma), P={_A\End_B}(\Sigma^\dagger)^{\rm op},
N={_B\Hom_B}(B,\Sigma\ot_A\Sigma^\dagger), M={_A\Hom_A}(\Sigma^\dagger\ot_B\Sigma,A),
\newtriangle, \newblacktriangle)$.
Then it follows from \thref{Morita-adjoint} that$(A,B,\Sigma^\dagger,\Sigma,\varepsilon,\eta)$ is a comatrix coring context if and only if 
$\MM(\Sigma,\Sigma^\dagger)$ is strict by a pair of invertible elements $(\varepsilon,\eta)$, formed by the counit of the corresponding comatrix coring, and the
unit of the corresponding matrix ring.

Assume more general that $\MM(\Sigma,\Sigma^\dagger)$ is strict.
Then there exist unique elements $\sum_{i\in I} \eta_i\ot_P\varepsilon_i\in N\ot_PM$ and
$\sum_{j\in J} \varepsilon'_j\ot_Q\eta'_j\in M\ot_QN$ such that
$\sum_i\eta_i\newtriangled \varepsilon_i=1_Q$ and $\sum_j \varepsilon'_j\newblacktriangled \eta'_j=1_P$.
Then $\cc=\Sigma^\dagger\ot_B\Sigma$ 
is no longer an $A$-coring; however, $\cc$ has local comultiplications and local counits
(compare to \cite{V:locunit}). For every $i\in I$ and $j\in J$, we define $\Delta_i=\Sigma^\dagger\ot_B\eta_i\ot_B\Sigma$ and 
$\Delta_j=\Sigma^\dagger\ot_B\eta'_j\ot_B\Sigma$. It is easily verified that
$(\cc\ot_A\Delta_k)\circ\Delta_l=(\Delta_l\ot_A\cc)\circ\Delta_k$ for all $k,l\in I\cup J$.
This means that the $\Delta_i$ and $\Delta_j$ are coassociating coassociative maps.
Moreover, they satisfy the generalized counit condition 
$$\sum_{i\in I}(\cc\ot_A\varepsilon_i)\circ\Delta_i=\cc=\sum_{j\in J}(\varepsilon_j\ot_A\cc)\circ\Delta_j=\cc.$$ 

In a similar way, $R=\Sigma\ot_A\Sigma^\dagger$ is a $B$-ring with local units and local multiplications.
The multiplications are defined as  $\mu_i=\Sigma\ot_A\varepsilon_i\ot_A\Sigma^\dagger$
and $\mu_j=\Sigma\ot_A\varepsilon'_j\ot_A\Sigma^\dagger$. Then we obtain $\mu_k\circ(R\ot_B\mu_l)=\mu_l\circ(\mu_k\ot_BR)$ for all $k,l\in I\cup J$. The generalized unit condition reads as
$$\sum_{i\in I}\mu_i\circ(R\ot_B\eta_i)=R=\sum_{j\in J}\mu_j\circ(\eta_j\ot_BR).$$

\subsection{Locally adjoint functors}\selabel{locadjfun}

We now consider $\CAT$, the bicategory whose $0$-cells are categories, $1$-cells are functors and $2$-cells are natural transformations. Then $\CAT$ is in fact even a $2$-category.
To avoid set-theoretical problems, we will consider a sub-bicategory $\Bb$ of $\CAT$, such that the natural transformations between each pair of functors form a set. In other words, $\Bb$, being enriched over $\Set$, fits into the setting of \seref{bicat}. 
Recall from the beginning of \seref{bicat} our convention to write the composition of $1$-cells and the horizontal composition of $2$-cells in a bicategory. This has important implications if we compute  the composition of functors and the horizontal composition of natural transformations in $\CAT$. Let $\Aa,\Bb$ and $\Cc$ be categories and $F:\Aa\to \Bb$ and $G:\Bb\to \Cc$ functors. Then we will denote 
\begin{equation}\eqlabel{covariant}
F\bullet_\Bb G= GF:\Aa\to\Cc
\end{equation}
for the composite functor. In the same way, for categories $\Aa,\Bb$ and $\Cc$, functors $F,G:\Aa\to \Bb$ and $H,K:\Bb\to \Cc$ and natural transformations $\alpha:F\to G$ and $\beta:H\to K$, we will denote
\begin{equation}\eqlabel{covariant2}
\alpha\ubul_\Bb\beta= \beta\alpha: HF\to KG,
\end{equation}
where the right hand side is the \emph{Godement product}\index{Godement
product} of natural transformations. 

 Take two categories $\Cc$ and $\Dd$, and two functors $F:\Cc\to \Dd$ and $G:\Dd\to \Cc$,
 and consider the Morita context from \thref{Morita}.
\begin{equation}\eqlabel{MorContextFunctors}
\MM(F,G)=(\Nat(G,G),\Nat(F,F)^{\rm op},\Nat(\Dd,FG),\Nat(GF,\Cc),\newdiamond,\newblackdiamond).
\end{equation}
The connecting maps are given by the following formulas, 
\begin{eqnarray*}
(\alpha\newdiamondd\beta)_D=\beta_{GD}\circ G\alpha_{D} &{\rm and}&
(\beta\newblackdiamondd\alpha)_C=F\beta_{C}\circ \alpha_{FC},
\end{eqnarray*}
where $\alpha\in\Nat(\Dd,FG)$, $\beta\in\Nat(GF,\Cc)$, $C\in\Cc$ and $D\in\Dd$.
By \thref{Morita-adjoint}, $(G,F)$ is an adjoint pair if and only if there exists a pair of invertible elements for the Morita context $\MM(G,F)$, i.e.\ if and only if we can find elements $\eta\in\Nat(\Dd,FG)$ and $\varepsilon\in\Nat(FG,\Cc)$ such that $\eta\newdiamondd\varepsilon=G$ and $\varepsilon\newblackdiamondd\eta=F$. 
Formulas (\ref{eq:triang1}-\ref{eq:triang2})
can be derived from this.

Applying left-right symmetry, we can construct a second Morita context, that describes the adjunction of the pair $(F,G)$:
\begin{equation}\eqlabel{MorContextFunctors2}
\bar{\MM}(F,G)=(\Nat(F,F),\Nat(G,G)^{\rm op},\Nat(\Cc,GF),\Nat(FG,\Dd),
\bar{\newdiamond},\bar{\newblackdiamond}),
\end{equation}
where
\begin{eqnarray*}
(\alpha\bar{\newdiamondd}\beta)_C=\beta_{FC}\circ F\alpha_{C} &{\rm and}&
(\beta\bar{\newblackdiamondd}\alpha)_D=G\beta_{D}\circ \alpha_{GD},
\end{eqnarray*}
for $\alpha\in\Nat(\Cc,GF)$ and $\beta\in\Nat(FG,\Dd)$.

We will now introduce the notion of a pair of locally adjoint functors.

\begin{Definitions}\delabel{locadjoint}
Consider functors $F: \Cc\to \Dd$, $G:\Dd\to \Cc$, and let
$\Ee$ be a generating subcategory for $\Cc$.

We call $G$ an \emph{$\Ee$-locally left adjoint} for $F$, if and only if, there exists a natural transformation $\varepsilon\in\Nat(GF,\Cc)$ and for all morphisms $f:E\to GD$ in $\Cc$, with $E\in \Ee$ and $D\in \Dd$, we can find a natural transformation $\eta^f\in\Nat(\Dd,FG)$, such that 
\begin{equation}
f=(\eta^f\newblackdiamondd\varepsilon)_D\circ f = \varepsilon_{GD}\circ G\eta^f_D\circ f.
\end{equation}
(In other words, the set $\Nat(\Dd,FG)\newblackdiamondd\varepsilon\subset \Nat(F,F)^{\rm op}$ acts with local units on $\Cc$.)

We call $F$ an \emph{$\Ee$-locally right adjoint} for $G$, if and only if, there exists a natural transformation $\varepsilon\in\Nat(GF,\Cc)$ and for all morphisms $f:FE\to FC$ in $\Dd$, with $E\in \Ee$ and $C\in \Cc$, we can find a natural transformation $\eta^f\in\Nat(\Dd,FG)$, such that 
\begin{equation}
f=(\varepsilon\newdiamondd\eta^f)_C\circ f = F\varepsilon_{C}\circ \eta^f_{FC}\circ f.
\end{equation}
(In other words, the set $\varepsilon\newdiamondd\Nat(\Dd,FG)\subset\Nat(G,G)$ acts with local units on $\Cc$.) 

If $F$ an $\Ee$-locally right adjoint for $G$ and $G$ is an $\Ee$-locally left adjoint for $F$, then we call $(G,F)$ an \emph{$\Ee$-locally adjoint pair}. If $(F,G)$ is an adjoint pair and $(G,F)$ is an $\Ee$-locally adjoint pair, then we call $(F,G)$ an \emph{$\Ee$-locally Frobenius pair}.
\end{Definitions}

\begin{Definition} (compare to \cite[Definitions 2.1 and 2.2]{Guo})
We use the same notation as in \deref{locadjoint}.
Suppose the category $\Dd$ has coproducts and consider the functor $S:\Dd\to \Dd,S(D)=D^{(I)}$, where $I$ is a fixed index set.

We call $G$ a \emph{left $\Ee$-locally quasi-adjoint} for $F$ if and only if $G$ is a left $\Ee$-locally adjoint for $SF$. We call $F$ a \emph{right $\Ee$-locally quasi-adjoint} for $G$ if and only if $F$ is a right $\Ee$-locally adjoint for $GS$.

We call $(G,F)$ an \emph{$\Ee$-locally quasi-adjoint pair} if and only if $G$ is a left $\Ee$-locally quasi-adjoint for $F$ and at the same time $F$ is a right $\Ee$-locally quasi-adjoint for $G$. We call $(F,G)$ a \emph{$\Ee$-locally quasi-Frobenius pair} if $(F,G)$ is an adjoint pair and $(G,F)$ an $\Ee$-locally quasi-adjoint pair.
\end{Definition}

\section{The induction functor}\selabel{induction}

\subsection{Adjunctions}\selabel{inductionadjoint}

Let $\cc$ be an $A$-bimodule. It is well-known (see e.g.\ \cite[18.28]{BrzWis:book}) that $\cc$ is an $A$-coring if and only if the functor $-\ot_A\cc : \Mm_A\to \Mm_A$ is a comonad. This comonad functor induces a functor
$$\Gg=\Gg^\cc: \Mm_A\to \Mm^\cc,\quad\Gg^\cc(N)=N\ot_A\cc;$$
where we denote $N\in \Mm_A$. 
The induction functor $\Gg^\cc$ has both a left adjoint $\Ff^\cc$ (the forgetful functor) and a right adjoint $\Hh^\cc$. These are given by
\begin{eqnarray*}
\Ff=\Ff^\cc:\Mm^\cc\to \Mm_A,&& \Ff^\cc(M)=M;\\
\Hh=\Hh^\cc:\Mm^\cc\to \Mm_A,&&\Hh^\cc(M)=\Hom^\cc(\cc,M);
\end{eqnarray*}
here we denote $M\in\Mm^\cc$. 

The unit and counit of these adjunctions are given by
\begin{eqnarray*}
\eta_M : M\to \Gg\Ff(M) = M\ot_A\cc,&& \eta_M(m)=m_{[0]}\ot_Am_{[1]}; \\
\epsilon_N : \Ff\Gg(N)=N\ot_A\cc \to N,&& \epsilon_N(n\ot_A c)=n\varepsilon_\cc(c);
\end{eqnarray*}
and
\begin{eqnarray}\eqlabel{kappalambda}
\lambda_N : N \to \Hh\Gg(N)=\Hom^\cc(\cc,N\ot_A\cc),&& \lambda_N(n)(c)=n\ot_A c;\nonumber\\
\kappa_M : \Gg\Hh(M)=\Hom^\cc(\cc,M)\ot_A\cc\to M,&& \kappa_M(f\ot_A c)=f(c);
\end{eqnarray}
for all $M\in \Mm^\cc$, $N\in\Mm_A$.

Recall that a functor is said to be Frobenius if it has a right  adjoint that is at the same time a left adjoint. Since adjoint functors are unique up to natural isomorphism, the study of the Frobenius property of the induction functor is related to the description of the sets ($k$-modules)
\begin{eqnarray*}
V=\Nat(\Ff,\Hh) &{\rm and}& W=\Nat(\Hh,\Ff).
\end{eqnarray*}

\begin{proposition}\prlabel{VW}
There exist isomorphisms of $k$-modules
\begin{eqnarray}
\Nat(\Gg\Ff,\id_{\Mm^\cc})\cong V=\Nat(\Ff,\Hh)\cong\Nat(\id_{\Mm^\cc},\Gg\Hh); \eqlabel{NatV}\\
\Nat(\id_{\Mm_A},\Ff\Gg)\cong W=\Nat(\Hh,\Ff)\cong\Nat(\Gg\Hh,\id_{\Mm_A}). \eqlabel{NatW}
\end{eqnarray}
\end{proposition}

\begin{proof}
The isomorphisms \equref{NatV} follow directly form the adjunctions $(\Ff,\Gg)$ and $(\Gg,\Hh)$ if we apply \equref{thetaadjoint}.

To prove \equref{NatW}, take any $\alpha\in\Nat(\Hh,\Ff)$, and define $\alpha'\in\Nat(\id_{\Mm_A},\Ff\Gg)$ as 
$$\alpha'_N=\alpha_{\Gg N}\circ\lambda_N.$$
Conversely, for any $\beta\in\Nat(\id_{\Mm_A},\Ff\Gg)$ we define $\beta'\in\Nat(\Hh,\Ff)$ by
$$\beta'_M=\Ff\kappa_M\circ\beta_{\Hh M}.$$
If we compute $\alpha''$, we find $\alpha''_M=\Ff\kappa_M\circ\alpha_{\Gg\Hh M}\circ\lambda_{\Hh M}$. By the naturality of $\alpha$, we know that $\Ff\kappa_M\circ\alpha_{\Gg\Hh M}=\alpha_M\circ \Hh\kappa_M$. 
Applying adjointness identity \equref{triang2} on the adjunction $(\Gg,\Hh)$, we obtain $\Hh\kappa_M\circ\lambda_{\Hh M}=\Hh M$. Combining both identities, we find that $\alpha''_M=\alpha_M$.

Similarly, we find $\beta''=\beta$, making use of \equref{triang1} on the adjunction $(\Gg,\Hh)$ and the naturality of $\beta$. Finally, $\Nat(\Hh,\Ff)\cong\Nat(\Gg\Hh,\id_{\Mm_A})$ follows in the same way from the adjunction $(\Ff,\Gg)$.
\end{proof}

\subsection{Description of sets of natural transformations}\selabel{nattranf}

To give a further description of $V$, let $\Rr$ be an additive subset of $\*c$ such that $A\Rr\subseteq \Rr A$ and consider the $k$-modules
\begin{eqnarray*}
V_1&=&{^\cc\Hom^\cc}(\cc\otimes\cc,\cc);\\
V_2&=&\{\theta\in{_A\Hom_A}(\cc\otimes_A\cc,A)~|~c_{(1)}\theta(c_{(2)}\otimes_Ad)=
\theta(c\ot_Ad_{(1)})d_{(2)}\};\\
V_3&=&{_A\Hom_{\*c}}(\cc,\*c);\\
V'_3&=&{_{\cc^*}\Hom_A}(\cc,\cc^*);\\
V_4&=&{_A\Hom^\cc}(\cc,\Rat_\Rr(\*c)),~\textrm{only if $\cc$ is $\Rr$-locally projective as a left $A$-module};\\
V_5&=&{_A\Hom^\cc(\cc,\*c)},~\textrm{only if $\cc$ is finitely generated and projective as a left $A$-module.}
\end{eqnarray*}

\begin{proposition}\prlabel{V3V'3}
Let $\cc$ be an arbitrary $A$-coring, then we have an isomorphism of $k$-modules
$$V_3={_A\Hom_{\*c}}(\cc,\*c)\cong V'_3={_{\cc^*}\Hom_A}(\cc,\cc^*).$$
\end{proposition}

\begin{proof}
Take any $\varphi\in {_A\Hom_{\*c}(\cc,\*c)}$, then we can easily construct a map
\[
\begin{array}{rrcl}
\widetilde{\varphi}:&\cc&\to &\cc^*;\\
&c&\mapsto & (d\mapsto \widetilde{\varphi}(c)(d)=\varphi(d)(c)).
\end{array}
\]
It is straightforward to check that switching the arguments as above corresponds in an isomorphism $V_3\cong V'_3$.
\end{proof}

By \cite[section 3.3]{CMZ} (see also \seref{coproduct} of this paper for a more general setting), $V\cong V_1\cong V_2$ for all corings, and $V\cong V_5$ if $\cc$ is finitely generated and projective as a left $A$-module.
We extend this result.

\begin{lemma}
\lelabel{formulecf}
Let $\cc$ be an $A$-coring which is $\Rr$-locally projective as a left $A$-module.
Then the following identity holds for all $f\in \Rat_\Rr(\*c)$ and $c\in\cc$~:
$$c_{(1)}f(c_{(2)})=f_{[0]}(c)f_{[1]}.$$
\end{lemma}

\begin{proof}
For all $g\in \Rr$, we have
\[
\begin{array}{rcl}
g(c_{(1)}f(c_{(2)})) & =&  (f*g)(c)\\
&=&(f_{[0]}g(f_{[1]}))(c)\\
&=&f_{[0]}(c)g(f_{[1]})\\
&=&g(f_{[0]}(c)f_{[1]}).
\end{array}
\]
Let now $\sum c_i\ot_Ag_i\in\cc\ot_A\Rr$ be a local basis for the elements $c_{(1)}f(c_{(2)})$ and $f_{[0]}(c)f_{[1]}$, then 
$c_{(1)}f(c_{(2)})=\sum c_ig_i(c_{(1)}f(c_{(2)}))=\sum c_ig_i(f_{[0]}(c)f_{[1]})=f_{[0]}(c)f_{[1]}$.
\end{proof}

\begin{proposition}
\prlabel{vcongv3}
Let $\cc$ be an arbitrary $A$-coring. Then there exists a map $\alpha_2:V\cong V_2\to V_3$.
If $\cc$ is $\Rr$-locally projective as a left $A$-module, then $V\cong V_3\cong V_4$. 
In particular, $V\cong V_5$ if $\cc$ is finitely generated and projective as a left $A$-module.
\end{proposition}

\begin{proof}
Let $\cc$ be any $A$-coring; then we can define a map
\[
\begin{array}{rrcl}
\alpha_2:&V_2&\to&V_3;\\
&\theta&\mapsto&(c\mapsto(d\mapsto\theta(d\otimes_Ac))).\\
\end{array}
\]
We verify that $\alpha_2$ is well-defined. First check $\alpha_2(\theta)=\bar{\varphi}$ is an $A$-bimodule map.
\[
\begin{array}{rcl}
(a\bar{\varphi}(c))(d)&=&\bar{\varphi}(c)(da)=\theta(da\ot_Ac)\\
&=&\theta(d\ot_Aac)=\bar{\varphi}(ac)(d);\\
(\bar{\varphi}(c)a)(d)&=&(\bar{\varphi}(c)(d))a=\theta(d\ot_Ac)a\\
&=&\theta(d\ot_Aca)=\bar{\varphi}(ca)(d).
\end{array}
\]
Next, we prove $\bar{\varphi}$ is also a right $\*c$-module map. Take $f\in \*c$, then we find
\[
\begin{array}{rcl}
\bar{\varphi}(c\cdot f)(d)&=& \bar{\varphi}(c_{(1)}\cdot f(c_{(2)}))(d)
=\bar{\varphi}(c_{(1)})(d)f(c_{(2)})\\
&=&f(\bar{\varphi}(c_{(1)})(d)c_{(2)})=f(\theta(d\ot_Ac_{(1)}c_{(2)}))\\
&=&f(d_{(1)}\theta(d_{(2)}\ot_Ac))=f(d_{(1)}\bar{\varphi}(c)(d_{(2)}))\\
&=&(\bar{\varphi}(c)*f)(d).
\end{array}
\]

Conversely, we can define a map
\[
\begin{array}{rrcl}
\alpha'_2:&V_3&\to&{_A\Hom_A}(\cc\ot_A\cc,A);\\
&\bar{\varphi}&\mapsto&(d\otimes_Ac\mapsto\varphi(c)(d)).
\end{array}
\]
We demonstrate that $\alpha'_2$ is well-defined. Take $\bar{\varphi}\in V_3$, then $\theta=\alpha'_2(\bar{\varphi})$ is an $A$-bimodule map:
\[
\begin{array}{rcl}
\theta(ad\otimes_Ac)&=&\bar{\varphi}(c)(ad)=a(\bar{\varphi}(c)d)\\
&=&a\theta(d\ot_Ac);\\
\theta(d\ot_Aca)&=&\bar{\varphi}(ca)(d)=(\bar{\varphi}(c)a)(d)\\
&=&(\bar{\varphi}(c)(d))a=\theta(d\ot_Ac)a.
\end{array}
\]
Now suppose that $\cc$ is $\Rr$-locally projective as a left $A$-module. We prove that the image of $\alpha_2$ lies within $V_2$.
Since $\bar{\varphi}$ is a $\*c$-module map, it follows from the theory of rational modules (see \cite{CVW2} and \cite{Wis:com}) that $\bar{\varphi}(d)\in\Rat_\Rr(\*c)$ and $\bar{\varphi}$ is also a $\cc$-comodule map between $\cc$ and $\Rat_\Rr(\*c)$. We can compute
\[
\begin{array}{rcl}
c_{(1)}\theta(c_{(2)}\ot_Ad)&=&c_{(1)}\bar{\varphi}(d)(c_{(2)})\\
&=&\bar{\varphi}(d)_{[0]}(c)\bar{\varphi}(d)_{[1]}\\
&=&\bar{\varphi}(d_{(1)})(c)d_{(2)}\\
&=&\theta(c\ot_Ad_{(1)})d_{(2)}.
\end{array}
\]
The second equation follows by \leref{formulecf} and the third one by the $\cc$-colinearity of $\bar{\varphi}$.

All the other implications are now straightforward.
\end{proof}

We will now describe the set $W$. Consider the following $k$-modules.
\begin{eqnarray*}
W_1&=&{_A\Hom_A}(A,\cc);\\
W_2&=&\{z\in\cc~|~az=za\}=\cc^A;\\
W_3&=&{_A\Hom_{\*c}}(\*c,\cc);\\
W'_3&=&{_{\cc^*}\Hom_A}(\cc^*,\cc);\\
W^r_3&=&{_A\Hom_{\*c}}(\Rat_\Rr(\*c),\cc),\textrm{~only if $\cc$ is $\Rr$-locally projective over $A$};\\
W_4&=&{_A\Hom^\cc}(\*c,\cc),\textrm{~only if $\cc$ is finitely generated and projective over $A$};\\
W_5&=&{_A\Hom_A}(\*c,A);\\
W^r_5&=&{_A\Hom_A}(\Rat_\Rr(\*c),A), \textrm{~only if $\cc$ is $\Rr$-locally projective over $A$}.
\end{eqnarray*}

Again by \cite[section 3.3]{CMZ} (or \seref{coproduct} of this paper), we know $W\cong W_1\cong W_2$ for arbitrary corings and $W\cong W_4$ if $\cc$ is finitely generated and projective as a left $A$-module. This can be easily generalized in the following way.

\begin{proposition}
\prlabel{wcongw3}
Let $\cc$ an $A$-coring, then $W\cong W_3\cong W'_3$. 
If $\cc$ is $\Rr$-locally projective as a left $A$-module, then $W^r_3\cong W^r_5$.
Consequently, if $\cc$ is finitely generated and projective as a left $A$-module, $(W\cong )W_3  \cong W_4\cong W_5$.
\end{proposition}

\begin{proof} By \leref{M^B} we immediately obtain that $W_3\cong W_2\cong W'_3$. Suppose that $\cc$ is $\Rr$-locally projective as a left $A$-module. By rationality properties we find that $W^r_3\cong{_A\Hom^\cc}(\Rat_R(\*c),\cc)$, and from the adjunction between the forgetful functor ${_A\Mm^\cc}\to{_A\Mm_A}$ and  $-\ot_A\cc$ we find that ${_A\Hom^\cc}(\Rat_R(\*c),\cc)\cong W^r_5$.
\end{proof}

To finish this section, we will describe the following classes of natural transformations
\begin{eqnarray}\eqlabel{XY}
X=\Nat(\Ff,\Ff) & {\rm and } &
Y=\Nat(\Gg,\Gg).
\end{eqnarray}

\begin{proposition}\prlabel{XY}
Let $\cc$ be an $A$-coring and consider the classes of natural transformations $X$ and $Y$ as in \equref{XY}. Then the following isomorphisms hold,
\begin{eqnarray*}
& X \cong \Nat(\Mm^\cc,\Gg\Ff) \cong \Nat(\Ff\Gg,\Mm_A) \cong Y \\
& \cong {^\cc\Hom^\cc}(\cc,\cc\ot_A\cc) \cong {_A\End^\cc}(\cc) \cong {^\cc\End_A}(\cc) \\
& \cong \*c^*
\cong(\*c)^A
\cong(\cc^*)^A\\
&\cong {_A\End_{\*c}}(\*c)\cong {_{\cc^*}\End_A}(\cc^*)\cong Z=\Nat(\Hh,\Hh).
\end{eqnarray*}
In particular, $X$, $Y$ and $Z$ are sets. 
\end{proposition}

\begin{proof}
The isomorphisms $X\cong \Nat(\Mm^\cc,\Gg\Ff)$ and $Y\cong\Nat(\Ff\Gg,\Mm_A)$ follow directly as an application of \equref{thetaadjoint} as $(\Ff,\Gg)$ is an adjoint pair. 

Take any $\alpha\in X$ and $(M,\rho^M)\in\Mm^\cc$. For any $m\in M$, the map 
$f_m:\cc\to M\ot_A\cc,\ c\mapsto m\ot_Ac$ is right $\cc$-colinear. We obtain by naturality of $\alpha$ that 
$\alpha_{M\ot_A\cc}=M\ot_A\alpha_\cc$. The naturality of $\alpha$ implies as well that $\rho^M\circ\alpha_M=\alpha_{M\ot_A\cc}\circ \rho^M$. This way we find
\begin{eqnarray*}
\alpha_M&=&(M\ot_A\varepsilon_\cc)\circ\rho^M\circ\alpha_M\\
&=&(M\ot_A\varepsilon_\cc)\circ\alpha_{M\ot_A\cc}\circ \rho^M\\
&=&(M\ot_A\varepsilon_\cc)\circ(M\ot_A\alpha_{\cc})\circ \rho^M.
\end{eqnarray*}
We conclude that $\alpha$ is completely determined by $\alpha_\cc$. By definition $\alpha_\cc\in\Hom_A(\cc,\cc)$ and by the naturality of $\alpha$ we find that $\alpha_\cc$ is left $\cc$-colinear as well. One can now easily see that the correspondence we obtained between $X$ and ${^\cc\End_A}(\cc)$ is bijective. 

Now take $\beta\in Y$. In a similar way as above, one can prove that $\beta_N=N\ot\beta_A$ for all $N\in\Mm_A$. Observe that by definition $\beta_A\in\End^\cc(\cc)$ and $\beta_A$ is left $A$-linear by the naturality of $\beta$. We conclude on the isomorphism $Y\cong{_A\End^\cc}(\cc)$.

The isomorphism $\End^\cc(\cc)\cong\*c$ restricts in a straightforward way to an isomorphism ${_A\End^\cc}(\cc)\cong\*c^*$ and similarly ${^\cc\End_A}(\cc)\cong\*c^*$.

Take $f\in(\*c)^A$, then for all $c\in\cc$, we find
$f(ca)=(af)(c)=(fa)(c)=f(c)a$, i.e.\ $f$ is right $A$ linear. This way we find that $\*c^*\cong(\*c)^A$ and dually $\*c^*\cong (\cc^*)^A$.

Furthermore, for any $\gamma\in{_A\End^\cc}(\cc)$, define $\theta\in{^\cc\End^\cc}(\cc,\cc\ot_A\cc)$ as $\theta(c)=c_{(1)}\ot_A\theta(c_{(2)})$ and conversely $\gamma=(\varepsilon\ot_A\cc)\circ\theta$.

Consider the map $\nu:(\*c)^A\to {_A\End_{\*c}}(\*c),\ \nu(f)(g)=f*g$, which has an inverse by evaluating at $\varepsilon$.

Finally, 
take $f\in {\*c^*}$. Then we define $\gamma\in Z$ as follows
$\gamma_M :\Hom^\cc(\cc,M)\to \Hom^\cc(\cc,M),\ \gamma_M(\varphi)(c)=\varphi(f\cdot c)$.
One easily checks that $\gamma_M$ is well-defined and natural in $M$. In this way we obtain a map $z:{^*\cc^*}\to Z$.
Conversely, for $\gamma\in Z$, take $\gamma_\cc(\cc)\in\Hom^\cc(\cc,\cc)$. Then by naturality of $\gamma$ one can easily check that $\gamma_\cc(\cc)$ is left $A$-linear. This way, we can define a map $z':Z\to {^*\cc^*},\ z'(\gamma)=\varepsilon_\cc\circ\gamma_\cc(\cc)$. Let us check that $z$ and $z'$ are each other inverses. For all $c\in{^*\cc^*}$,
$z'\circ z(f)(c)=\varepsilon(f\cdot c)=f(c)$. For all $\gamma\in Z$, $M\in\Mm^\cc$ and $\varphi\in\Hom^\cc(\cc,M)$ we find
\begin{eqnarray*}
z\circ z'(\gamma_M)(\varphi)(c)&=&\varphi((\varepsilon\circ\gamma_\cc(\cc))\cdot c)\\
&=&\varphi(\gamma_\cc(\cc)(c))=(\gamma_M(\varphi))(c),
\end{eqnarray*}
where the last equation follows from the naturality of $\gamma$, applied to the morphism $\varphi\in\Hom^\cc(\cc,M)$.
\end{proof}

\begin{remark}\relabel{XY}
The above theorem only states isomorphisms of modules. However, some of these objects have an additional ring structure. All stated (iso)morphisms are also ring morphisms for those objects that posses a ring structure, but sometimes one has to consider the opposite multiplication.
For sake of completeness, we state the correct isomorphisms, but we leave the proof to the reader.
\begin{eqnarray*}
&\Nat(\Ff,\Ff)^{\rm op}\cong \Nat(\Gg,\Gg)\cong \Nat(\Hh,\Hh)^{\rm op}\\
&\cong{^\cc\End_A(\cc)}^{\rm op} \cong {_A\End^\cc}(\cc)\cong {^*\cc^*}\\
&\cong {_A\End_{\*c}}(\*c)\cong {_{\cc^*}\End_A}(\cc^*)^{\rm op}\\
\end{eqnarray*}
\end{remark}

\subsection{The Yoneda-approach}\selabel{yoneda}

\begin{lemma}\lelabel{contravariant}
\begin{enumerate}
\item Let $N$ be a $\cc^*\hbox{-}A$ bimodule. Then the following assertions hold
\begin{enumerate}
\item ${_{\cc^*}\Hom_A}(N,\cc^*)\in{\Mm_{\*c^*}}$;
\item ${_{\cc^*}\Hom_A}(N,\cc)\in{\Mm_{\*c^*}}$.
\end{enumerate}
\item Let $M$ be a $\cc\hbox{-}A$ bicomodule and $\Rr\subset\cc^*$. Then the following assertions hold
\begin{enumerate}
\item ${^\cc\Hom_A}(M,\cc)\cong {_A\Hom_A}(M,A)\in{\Mm_{\*c^*}}$;
\item ${_{\cc^*}\Hom_A}(M,\cc^*)\cong{_{\cc^*}\Hom_A}(M,{_\Rr\Rat}(\cc^*))\cong{^\cc\Hom_A}(M,{_\Rr\Rat}(\cc^*)) \in{\Mm_{\*c^*}}$;
\end{enumerate}
\end{enumerate}
where the $\Rr$-rational part of $\cc^*$ is only considered if $\cc$ is $\Rr$-locally projective as a right $A$-module. 
\end{lemma}

\begin{proof}
$\ul{(1a)}$ Take any $N\in{_{\cc^*}\Mm_A}$, for any $f\in{_{\cc^*}\Hom_A}(N,\cc^*)$, $g\in\*c^*$ and $n\in N$, we define
\begin{equation}\eqlabel{action}
(f*g)(n)=(f(n))*g.
\end{equation}
Note that $f*g$ is right $A$-linear, since $g$ commutes with all elements of $A$ by \prref{XY}. One can easily verify that \equref{action} defines a left $\*c^*$-action.

$\ul{(1b)}$ We give only the explicit form of the action and leave other verifications to the reader.
Take any $f\in{_{\cc^*}\Hom_A}(N,\cc)$, $g\in\*c^*$ and $n\in N$, then we define
\begin{equation}\eqlabel{action2}
(f*g)(n)=f(n)\cdot g=f(n)_{(1)} g(f(n)_{(2)}).
\end{equation}

$\ul{(2a)}$
Analogously to the adjunction of $(\Ff^\cc,\Gg^\cc)$, the forgetful functor ${^\cc\Mm_A}\to {_A\Mm_A}$ has a right  adjoint $-\ot_A\cc:{_A\Mm_A}\to {^\cc\Mm_A}$. Consequently, for any $M\in{^\cc\Mm_A}$, we have an isomorphism ${^\cc\Hom_A}(M,\cc)\cong {_A\Hom_A}(M,A)$ that is natural in $M$. 
Moreover, the action defined in \equref{action2} can be restricted to a right $\*c^*$-module structure on  ${^\cc\Hom_A}(M,\cc)\cong{_A\Hom_A}(M,A)$. For any $f\in{_A\Hom_A}(M,A)$, $m\in M$ and $g\in\*c^*$, one defines explicitly 
$$(f*g)(m)= g(m_{[-1]}f(m_{[0]}))=g(m_{[-1]})f(m_{[0]}).$$

$\ul{(2b)}$
Since every left $\cc$-comodule is also a left $\*c$-module (see \seref{modversuscomod}), by part $(1a)$, we find that ${_{\cc^*}\Hom_A}(M,\*c)\in{\Mm_{\*c^*}}$.

Suppose now that $\cc$ is $\Rr$-locally projective as a right  $A$-module. Then the image of any $f\in{_{\cc^*}\Hom_A}(M,\*c)$ lies within the rational part ${_\Rr\Rat}(\cc^*)$. Indeed, for any $g\in \cc^*$, $g* f(m)= f(g\cdot m)=f(g(m_{[-1]})m_{[0]})=g(m_{[-1]})f(m_{[0]})$, so $f(m)\in {_\Rr\Rat}(\cc^*)$. 
We can conclude that ${_{\cc^*}\Hom_A}(M,\cc^*)\cong{_{\cc^*}\Hom_A}(M,{_\Rr\Rat}(\cc^*))\cong{^\cc\Hom_A}(M,{_\Rr\Rat}(\cc^*))$.
\end{proof}

The observations made in \leref{contravariant} lead to the introduction of the following contravariant functors
\begin{eqnarray}
\Jj :&{_{\cc^*}\Mm_A} \to {\Mm_{\*c^*}}, & \Jj(M)={_{\cc^*}\Hom_A}(M,\cc^*);\nonumber\\
\Kk :&{_{\cc^*}\Mm_A} \to {\Mm_{\*c^*}}, & \Kk(M)={_{\cc^*}\Hom_A}(M,\cc);\nonumber\\
\Jj' :& {^\cc\Mm_A} \to {\Mm_{\*c^*}},& \Jj'(M)={_{\cc^*}\Hom_A}(M,\cc^*)\cong{_{\cc^*}\Hom_A}(M,{_\Rr\Rat}(\cc^*))\nonumber\\
&&\hspace{2cm}\cong{^\cc\Hom_A}(M,{_\Rr\Rat}(\cc^*));\nonumber\\
\Kk' :& {^\cc\Mm_A} \to {\Mm_{\*c^*}},&\Kk'(M)= {_A\Hom_A}(M,A)\cong{^\cc\Hom_A}(M,\cc).\eqlabel{JK}
\end{eqnarray}
(The alternative descriptions of $\Jj'$ in terms of the $\Rr$-rational part of $\cc^*$ is only considered if $\cc$ is $\Rr$-locally projective as a right $A$-module.)
Out of these functors we can construct the $k$-modules
\begin{eqnarray*}
V_6=\Nat(\Kk',\Jj') &{\rm and}&
V_7=\Nat(\Kk,\Jj);\\
W_6=\Nat(\Jj',\Kk') &{\rm and}&W_7=\Nat(\Jj,\Kk).
\end{eqnarray*}

\begin{lemma}\lelabel{covariant}
Let $N$ be an $\cc^*\hbox{-}A$ bimodule. Then
\begin{enumerate}
\item 
$ N^A
\cong{_A\Hom_A}(A,N)
\cong {_{\cc^*}\Hom_A}(\cc^*,N)\in {_{\*c^*}\Mm}$;
\item 
${_{\cc^*}\Hom_A}(\cc,N)\in{_{\*c^*}\Mm}$.
\end{enumerate}
\end{lemma}
\begin{proof}
$\ul{(1)}$ 
Both isomorphisms follow from \leref{M^B}.
We define a left $\*c^*$-action on ${_{\cc^*}\Hom_A}(\cc^*,N)$ with the following formula
$$(f*\varphi)(g)=\varphi(g*f)$$
for all $\varphi\in{_{\cc^*}\Hom_A}(\cc^*,N)$, $f\in \*c^*$ and $g\in\cc^*$.

$\ul{(2)}$ For $\varphi\in{_{\cc^*}\Hom_A}(\cc,N)$, $f\in \*c^*$ and $c\in\cc$ we define
$$(f*\varphi)(c)=\varphi(c_{(1)}f(c_{(2)})).$$
One can easily verify this turns ${_{\cc^*}\Hom_A}(\cc,N)$ into a left $\*c^*$-module.
\end{proof}

Using \leref{covariant} we can construct the covariant functors
\begin{eqnarray}
\widetilde{\Jj} : &{_{\cc^*}\Mm_A} \to {_{\*c^*}\Mm}, &\widetilde{\Jj}(N)=
N^A\cong{_A\Hom_A}(A,N)\cong 
{_{\cc^*}\Hom_A}(\cc^*,N);\nonumber\\
\widetilde{\Kk} :& {_{\cc^*}\Mm_A} \to {_{\*c^*}\Mm};&\widetilde{\Kk}(N)={_{\cc^*}\Hom_A}(\cc,N)\eqlabel{JKtilde}
\end{eqnarray}
and the $k$-modules
\begin{eqnarray*}
V_8=\Nat(\widetilde{\Jj},\widetilde{\Kk})&{\rm and}&W_8=\Nat(\widetilde{\Kk},\widetilde{\Jj}).
\end{eqnarray*}

Let $\Xx$ be any category, $F:\Xx\to \Set$ a covariant functor and $X\in\Xx$. Recall that by the Yoneda Lemma (see e.g.\ \cite[Theorem 1.3.3]{Bor1}) $\Nat(\Hom(X,-),F)\cong F(X)$. Similarly for any contravariant functor $G:\Xx\to \Set$, we have $\Nat(\Hom(-,X),G)\cong G(X)$. Of course the Yoneda Lemma can be applied to the particular case where $F=\Hom(X,-)$ and $G=\Hom(-,X)$. In those cases, $\Nat(F,F)$ and $\Nat(G,G)$ can be completed with a semigroup structure, coming from the composition of natural transformations. The following Lemma compares these structures with the semigroup structure of $\Hom(X,X)$ (under composition). This result might be well-known, but since we could not find any reference, we include the proof.

\begin{lemma}\lelabel{Yonedaring}
Let $\Xx$ be any category and $X\in\Xx$, then we have the following isomorphisms of semigroups
$$\Nat(\Hom(X,-),\Hom(X,-))^{\rm op}\cong\Hom(X,X)\cong\Nat(\Hom(-,X),\Hom(-,X)).$$
\end{lemma}

\begin{proof}
Consider the Yoneda bijection $\Lambda:\Nat(\Hom(X,-),\Hom(X,-))\to \Hom(X,X)$;  $\Lambda(\alpha)=\alpha_X(X)$. Let us compute $\Lambda(\alpha\circ\beta)=(\alpha\circ\beta)_X(X)=\alpha_X\circ\beta_X(X)$. Consider the morphism $\beta_X(X):X\to X$ and apply the naturality of the functor $\Hom(X,-)$ to this morphism, we obtain $\alpha_X\circ\beta_X(X)=\beta_X(X)\circ\alpha_X(X)$.

Similarly, starting from the bijection ${\rm V}:\Nat(\Hom(-,X),\Hom(-,X))\to \Hom(X,X)$;  ${\rm V}(\alpha)=\alpha_X(X)$, we find
${\rm V}(\alpha\circ\beta)=(\alpha\circ\beta)_X(X)=\alpha_X\circ\beta_X(X)$. The functor property of the contravariant functor $\Hom(-,X)$ implies $\alpha_X\circ\beta_X(X)=\alpha_X(X)\circ\beta_X(X)$ and we find the needed semigroup morphism.
\end{proof}

\begin{proposition}\prlabel{vw6}
Let $\cc$ be an $A$-coring. Then we have isomorphisms of $k$-modules
\begin{enumerate}
\item $V'_3\cong V_7\cong V_8$;
\item $W'_3\cong W_7\cong W_8$;
\item $\Nat(\Jj,\Jj)\cong \Nat(\widetilde{\Jj},\widetilde{\Jj})^{\rm op}\cong {_{\cc^*}\End_A}(\cc^*)\cong X $;
\item $\Nat(\Kk,\Kk)\cong \Nat(\widetilde{\Kk},\widetilde{\Kk})^{\rm op}\cong {_{\cc^*}\End_A}(\cc)$;
\item $\Nat(\Kk',\Kk')\cong {_{\cc^*}\End_A}(\cc)$.
\end{enumerate}
If $\cc$ is $\Rr$-locally projective as a left $A$-module, then 
$V_6 \cong V'_3$, 
$W_6\cong {_A\Hom^\cc}(\cc,\Rat_\Rr(\*c))$, $\Nat(\Jj',\Jj') \cong {_A\End^\cc}(\Rat_\Rr(\*c))$ and $\Nat(\Kk',\Kk')\cong\Nat(\Kk,\Kk)$.
\end{proposition}

\begin{proof}
All isomorphisms are immediate consequences of the Yoneda Lemma and Lemma \ref{le:Yonedaring}.
\end{proof}

\subsection{The coproduct functor}\selabel{coproduct}

Quasi-Frobenius type properties can not be described by the functors $\Ff$ and $\Gg$ alone, we have to incorporate a new functor in our theory (compare also with \cite{Guo}). 

Consider the following coproduct-functor
$$\Ss:\Mm_A\to \Mm_A,\ \Ss(M)=M^{(I)},$$
where $I$ is an arbitrary fixed index set.

Applying our previous results, we will give a full description of the sets 
\begin{eqnarray*}
\Nat(\Ss\Ff,\Ss\Ff),
&\Nat(\Mm_A,\Ss\Ff\Gg),
&\Nat(\Gg\Ss\Ff,\Mm^\cc).
\end{eqnarray*}

To improve the readability of the next theorems, let us recall the construction of coproducts in $\Mm^{\cc}$. Take $(M,\rho_M)\in\Mm^\cc$; then $(M,\rho_M)^{(I)}=(M^{(I)},\rho)$, where the coaction $\rho$ is given by the following composition
\begin{equation}\eqlabel{coactioncoproduct}
\xymatrix{
\rho: 
\Ss(M) \ar[rr]^{\Ss(\rho_M)} && 
\Ss(M\ot_A\cc) \ar[rr]^{\cong} && 
\Ss(M)\ot_A\cc
}
\end{equation}
where we used that the tensor product commutes with coproducts.

\begin{lemma}\lelabel{coproductHoms}
Let $\cc$ be an $A$-coring. Then we have the following isomorphisms of $k$-modules
$${_A\Hom^\cc}(\cc,\Ss(\cc))\cong{^\cc\Hom^\cc}(\cc,\Ss(\cc\ot_A\cc))
\cong{^\cc\Hom_A}(\cc,\Ss(\cc)).$$
\end{lemma}

\begin{proof}
Take $\gamma\in{_A\Hom^\cc}(\cc,\Ss(\cc))$. Then we define
\[
\xymatrix{
\theta: \cc \ar[rr]^-\Delta && \cc\ot_A\cc \ar[rr]^-{\cc\ot_A\gamma} && \cc\ot_A\Ss(\cc)\cong \Ss(\cc\ot_A\cc).
}\]
Conversely, given $\theta\in{^\cc\Hom^\cc}(\cc,\Ss(\cc\ot_A\cc))$, define
\[
\xymatrix{
\gamma: \cc \ar[rr]^-\theta && \Ss(\cc\ot_A\cc)\cong\cc\ot_A\Ss(\cc) \ar[rr]^-{\varepsilon\ot_A\cc} && \Ss(\cc).
}\]
The second isomorphism is constructed in the same way.
\end{proof}

\begin{proposition}\prlabel{SFSF}
There exist maps
\[
\xymatrix{
\Nat(\Ss\Ff,\Ss\Ff) \ar[r]^\cong & {^{\cc}\End_{A}}({\cc}^{(I)}) \ar[r]^{\upsilon} \ar[d]^{\cong} & {_{\cc^*}\End_{A}}({\cc}^{(I)})\\
\Nat(\Gg\Ss,\Gg\Ss)^{\rm op} \ar[r]^\cong
& {_A\End^{\cc}}({\cc}^{(I)})^{\rm op} \ar[r]^{\upsilon'} & {_A\End_{\*c}}({\cc}^{(I)})^{\rm op}
}
\]
where $\upsilon$ (resp. $\upsilon'$) is an isomorphism as well if $\cc$ is locally projective as a left (resp. right) $A$-module.
\end{proposition}

\begin{proof}
Take $\alpha\in\Nat(\Ss\Ff,\Ss\Ff)$. Then we find, by definition, that $\alpha_\cc\in \End_A(\Ss(\cc))$. Take now $N\in\Mm_A$. For any $n\in N$, we can consider the right $\cc$-colinear map $f_n:\cc\to N\ot_A\cc$, $f_n(c)=n\ot_Ac$. The naturality of $\alpha$ and the commutativity of the tensor product and coproduct imply the commutativity of the following diagram.
\[
\xymatrix{
N\ot_A\Ss(\cc)\ar[r]^{\cong} & \Ss(N\ot_A\cc) \ar[rr]^-{\alpha_{N\ot_A\cc}}&& \Ss(N\ot_A\cc) \ar[r]^{\cong} & N\ot_A\Ss(\cc)\\
& \Ss(\cc) \ar[rr]_{\alpha_\cc} \ar[ul]^{f_n\ot_A\Ss(\cc)} \ar[u]^{\Ss(f_n)} && \Ss(\cc) \ar[ur]_{f_n\ot_A\Ss(\cc)} \ar[u]_{\Ss(f_n)}
}
\]
This implies that $\alpha_{N\ot_A\cc}$ is determined by $\alpha_\cc$ up to isomorphism, as expressed in the following diagram.
\[\xymatrix{
\Ss(N\ot_A\cc) \ar[rr]^{\alpha_{N\ot_A\cc}} \ar[d]_{\cong} && \Ss(N\ot_A\cc) \\
N\ot_A\Ss(\cc) \ar[rr]_{N\ot_A\alpha_\cc} && N\ot_A\Ss(\cc) \ar[u]_{\cong}
}\]
It follows now easily from the naturality of $\alpha$ that $\alpha_\cc$ is left $\cc$-colinear, and thus $\alpha_\cc\in{^\cc\End_A}(\cc)$. Moreover, $\alpha$ is completely determined by its value in $\cc$. Take any $M\in\Mm^\cc$ and consider the following diagram.
\[
\xymatrix{
\Ss(M)\ar[rr]^-{\alpha_M} \ar[d]_{\Ss(\rho_M)} && \Ss(M) \ar[d]^{\Ss(\rho_M)}\\
\Ss(M\ot_A\cc)\ar[d]_{\cong} \ar[rr]^{\alpha_{M\ot_A\cc}}  && \Ss(M\ot_A\cc) \ar[d]^{\cong}\\
M\ot_A\Ss(\cc) \ar[rr]_{M\ot_A\alpha_\cc} && M\ot_A\Ss(\cc) \ar[r]_{\cong} & \Ss(M)\ot_A\cc \ar@(ur,r)[uul]_-{\Ss(M)\ot_A\varepsilon} 
}
\]
The upper quadrangle commutes by the naturality of $\alpha$, applied on the $\cc$-colinear morphism $\rho_M: M\to M\ot_A\cc$, the lower quadrangle commutes by the previous observations and the commutativity of the triangle is exactly the counit condition on the comodule $\Ss(M)$.
This way we find an isomorphism $\Nat(\Ss\Ff,\Ss\Ff)\cong{^\cc\End_A}(\cc^{(I)})$. The second horizontal isomorphism is proved in the same way.
The vertical isomorphism is a consequence of \leref{sumproduct} and \leref{coproductHoms}:
$${_A\Hom^\cc}(\cc^{(I)},\cc^{(I)})\cong({_A\Hom^\cc}(\cc,\cc^{(I)}))^I 
\cong ({^\cc\Hom_A}(\cc,\cc^{(I)}))^I\cong {^\cc\Hom_A}(\cc^{(I)},\cc^{(I)}).$$

We leave it to the reader that the constructed isomorphisms are algebra morphisms.
The morphisms $\upsilon$ and $\upsilon'$ follow from the relations between left $\cc$-comodules 
and left $\cc^*$-modules (see \seref{modversuscomod}).
\end{proof}

\begin{lemma}\lelabel{modiso's}
Let $\cc$ be an $A$-coring, $B\to A$ a ring morphism and $I$ any index set. 
\begin{enumerate}[(i)]
\item ${_BV_I}:={_B\Hom_{\*c}}(\cc,(\*c)^I)\cong {_B\Hom_{\*c}}(\cc^{(I)},\*c)\\
\texttt{ }~\textrm{ }~\textrm{ }\stackrel{\xi}{\cong} {_{\cc^*}\Hom_B}(\cc^{(I)},\cc^*)\cong{_{\cc^*}\Hom_B}(\cc,(\cc^*)^I)$; 
\item for all $M\in{^\cc\Mm^\cc}$, 
$${^\cc\Hom^\cc}(M,\cc)\cong\{\theta\in{_A\Hom_A}(M,A)~|~x_{[-1]}\theta(x_{[0]})=\theta(x_{[0]})x_{[1]}, \textrm{for all}\ x\in M\};$$
\item there exist morphisms $
{^\cc\Hom^\cc}(\Ss(\cc\ot_A\cc),\cc)\stackrel{\xi_1}\to {_AV_I}\cong V_3^I$, where $\xi_1$ becomes an isomorphism if $\cc$ is locally projective as a left $A$-module;
\item ${_B\Hom_{\*c}}((\*c)^{(I)},\cc)
{\cong}{_{\cc^*}\Hom_B}(\cc^*,\cc^{I})
{\cong}(\cc^I)^B\cong(\cc^B)^{I}$;
\item ${_BW_I}:={_{\cc^*}\Hom_B}(\cc^*,\cc^{(I)})\cong{_B\Hom_{\*c}}(\*c,\cc^{(I)})\cong (\cc^{(I)})^B\cong (\cc^B)^{(I)}$.
\end{enumerate}
\end{lemma}

\begin{proof} \ul{(i)} The first and last isomorphism are an immediate consequence of \leref{sumproduct}, the second isomorphism is induced by the isomorphism of \prref{V3V'3}.\\
\ul{(ii)} Take any $\gamma\in {^\cc\Hom^\cc}(M,\cc)$ and define $\theta=\varepsilon\circ\gamma$. Clearly, $\theta\in{_A\Hom_A}(M,A)$. Moreover, by the bi-colinearity of $\gamma$ we find for all $x\in M$,
$$x_{[1]}\ot_A\gamma(x_{[0]}) = \gamma(x)_{(1)}\ot_A \gamma(x)_{(2)} 
=\gamma(x_{[0]})\ot_A x_{[1]}.$$
If we apply $\cc\ot_A\varepsilon$ to the first equation, $\varepsilon\ot_A\cc$ to the second equation we obtain $x_{[-1]}\theta(x_{[0]})=\theta(x_{[0]})x_{[1]}=\gamma(x)$. Conversely, starting from $\theta\in {_A\Hom_A}(M,A)$, such that $x_{[-1]}\theta(x_{[0]})=\theta(x_{[0]})x_{[1]}$ for all $x\in M$, we define $\gamma(x)=x_{[-1]}\theta(x_{[0]})$.\\ 
\ul{(iii)} Denote by $\iota^{\cc\ot_A\cc}_\ell :\cc\ot_A\cc \to \Ss(\cc\ot_A\cc)$ and $\iota_\ell^\cc: \cc\to \Ss(\cc)$ the canonical injections. Consider the following diagrams.
\[
\xymatrix{
\Ss(\cc\ot_A\cc) \ar[rr]^{\nu} && \cc  & && \*c\\
\cc\ot_A\cc \ar[u]^{\iota^{\cc\ot_A\cc}_\ell} \ar[urr]_{\nu_\ell}
&& & \cc \ar[urr]^{\varphi_\ell} \ar[rr]_{\iota_\ell^\cc} && \Ss(\cc) \ar[u]_{\varphi}
}
\]
We find that every morphism $\nu\in {^\cc\Hom^\cc}(\Ss(\cc\ot_A\cc),\cc)$ is completely determined by the morphisms $\nu_\ell\in {^\cc\Hom^\cc}(\cc\ot_A\cc,\cc)=V_2$. Similarly, any $\varphi\in V_I^A$ is completely determined by $\varphi_\ell\in V_3$. The needed morphism and isomorphism is now a consequence of \prref{vcongv3} and \leref{sumproduct}.\\
\ul{(iv)} The first isomorphism is a consequence of \leref{sumproduct} the second one is a consequence of \leref{M^B}. The last isomorphism is trivial.\\
\ul{(v)} The second (and first) isomorphism follows from \leref{M^B}, the last one is trivial.
\end{proof}

\begin{proposition}\prlabel{SFG,GSF}
Let $\cc$ be an $A$-coring and $I$ any index set, than the following isomorphisms hold.
\begin{enumerate}[(i)]
\item $\Nat(\Mm_A,\Ss\Ff\Gg){\cong}(\cc^A)^{(I)}$;
\item $\Nat(\Gg\Ss\Ff,\Mm^\cc)\cong {^\cc\Hom^\cc}(\Ss(\cc\ot_A\cc),\cc)$.
\end{enumerate}
\end{proposition}

\begin{proof}
\ul{(i)}
Take any $\zeta\in\Nat(\Mm_A,\Ss\Ff\Gg)$. Then $\zeta_A\in{\Hom_A}(A,\cc^{(I)})$ by definition, and from the naturality of $\zeta$ we obtain that $\zeta_A$ is left $A$-linear. Applying the same techniques as in the proof of \prref{SFSF}, we find that $\zeta$ is completely determined by $\zeta_A$, and thus we obtain an isomorphism $\Nat(\Mm_A,\Ss\Ff\Gg){\cong}{_A\Hom_A}(A,\cc^{(I)})\cong (\cc^A)^{(I)}$.\\
\ul{(ii)} 
The proof is completely similar to part $(i)$. Any $\nu \in \Nat(\Gg\Ss\Ff)$ is completely determined by $\nu_\cc : \Ss\Ff\Gg(\cc)=\Ss(\cc\ot_A\cc)\to \cc$, by definition $\nu_\cc$ is right $\cc$-colinear and the left $\cc$-colinearity follows from the naturality of $\nu$, i.e.\ $\nu_\cc\in {^\cc\Hom^\cc}(\Ss(\cc\ot_A\cc),\cc)$.
\end{proof}

We give a generalization of \prref{VW}, the proof is completely similar.

\begin{proposition}\prlabel{VWS}
There exist isomorphisms of $k$-modules
\begin{enumerate}[(i)]
\item $\Nat(\Gg\Ss\Ff,\id_{\Mm^\cc})\cong\Nat(\Ss\Ff,\Hh);$
\item $\Nat(\id_{\Mm_A},\Ss\Ff\Gg)\cong\Nat(\Hh,\Ss\Ff).$
\end{enumerate}
\end{proposition}

\begin{proof}
\ul{(i)} The isomorphism follows directly form the adjunction between $\Gg$ and $\Hh$ if we apply \equref{thetaadjoint}.

\ul{(ii)}. Take any $\alpha\in\Nat(\Hh,\Ss\Ff)$, then we define $\alpha'\in\Nat(\id_{\Mm_A},\Ss\Ff\Gg)$ as 
$$\alpha'_N=\alpha_{\Gg N}\circ\lambda_N.$$
Conversely, for any $\beta\in\Nat(\id_{\Mm_A},\Ss\Ff\Gg)$ we define $\beta'\in\Nat(\Hh,\Ss\Ff)$ by
$$\beta'_M=\Ss\Ff\kappa_M\circ\beta_{\Hh M}.$$
If we compute $\alpha''$, we find 
\begin{eqnarray*}
\alpha''_M&=&\Ss\Ff\kappa_M\circ\alpha_{\Gg\Hh M}\circ\lambda_{\Hh}\\
&=&\alpha_M\circ \Hh\kappa_M\circ\lambda_{\Hh}=\alpha_M
\end{eqnarray*}
where we used the naturality of $\alpha$ in the second equality and \equref{triang2} on the adjunction $(\Gg,\Hh)$ in the third equality.
Similarly, we find
\begin{eqnarray*}
\beta''_M&=&\Ss\Ff\kappa_{\Gg N}\circ\beta_{\Hh\Gg N}\circ\lambda_{N}\\
&=&\Ss\Ff\kappa_{\Gg N}\circ\Ss\Ff\Gg\lambda_N\circ\beta_N\\
&=&\Ss\Ff(\kappa_{\Gg N}\circ\Gg\lambda_N)\circ\beta_N = \beta_N.
\end{eqnarray*}
Here we used the naturality of $\beta$ in the second equality and \equref{triang1} in the fourth equality.
\end{proof}

Consider the functors
\begin{eqnarray}\eqlabel{KkS}
\Kk^s :&{_{\cc^*}\Mm_A} \to {\Mm_{\*c^*}}, & \Kk^s(M)={_{\cc^*}\Hom_A}(M,\Ss(\cc));\nonumber\\
\widetilde{\Kk}^s : &{_{\cc^*}\Mm_A} \to {_{\*c^*}\Mm}, &\widetilde{\Kk}^s(N)=
{_{\cc^*}\Hom_A}(\Ss(\cc),N).
\end{eqnarray}
As a consequence of the Yoneda Lemma, we immediately obtain the following

\begin{proposition}\prlabel{NatYonS}
With notation as introduced before, the following isomorphisms hold:
\begin{enumerate}[(i)]
\item $\Nat(\Kk^s,\Kk^s)\cong \Nat(\widetilde{\Kk}^s,\widetilde{\Kk}^s)^{\rm op}\cong{_{\cc^*}\End_A}(\Ss(\cc))$;
\item $\Nat(\Kk^s,\Jj)\cong \Nat(\widetilde{\Jj},\widetilde{\Kk}^s)\cong {_{\cc^*}\Hom_A}(\Ss(\cc),\cc^*)$;
\item $\Nat(\widetilde{\Kk}^s,\widetilde{\Jj})\cong \Nat(\Jj,\Kk^s)\cong
{_{\cc^*}\Hom_A}(\cc^*,\Ss(\cc))$.
\end{enumerate}
\end{proposition}

\section{Characterizations of co-Frobenius and quasi-co-Frobenius corings}\selabel{Morita}

\subsection{Locally Frobenius corings}\selabel{frob}

\begin{lemma}\lelabel{imj}
Let $\cc$ be an $A$-coring and $B\to A$ a ring morphism. And take any $j\in{_B\Hom_{\*c}}(\cc,\*c)$. Then 
$\im j$ is a right ideal in $\*c$ and $(\im j)^B$ is a right ideal in $(\*c)^B$. 
\end{lemma}

\begin{proof}
Take $f\in \im j$, i.e.\ $f=j(c)$ for some $c\in\cc$. Then for any $g\in\*c$, $f*g=j(c)*g=j(c\cdot g)\in\im j$ by the right $\*c$-linearity of $j$.\\
Suppose now that $f=j(c)\in(\im j)^B$ and $g\in(\*c)^B$. We have to check that $j(c\cdot g)$ commutes with all $b\in B$. We find $bj(c\cdot g)=bj(c)* g=j(c)b *g=j(c)*gb=j(c\cdot g)b$.
\end{proof}

\begin{lemma}\lelabel{j'}
Let $\cc$ be an $A$-coring and $B\to A$ a ring morphism.
Consider $j\in{_B\Hom_{\*c}}(\cc,\*c)$. The restriction of $j$ on $\cc^B$ defines a map
$$j':\cc^B\to (\*c)^B$$
that is $Z(B)\hbox{-}(\*c)^B$-bilinear, where $Z(B)$ denotes the center of $B$. Moreover, $\im j'$ is a right ideal in $(\*c)^B$.
\end{lemma}

\begin{proof}
Take $c\in\cc^B$ and $b\in B$. Then $bj(c)=j(bc)=j(cb)=j(c)b$, so $j(c)\in(\*c)^B$. Since $j$ is $B\hbox{-}\*c$-bilinear and $\cc^A\subset\cc$ is a bimodule with restricted actions of $Z(B)\subset B$ and $(\*c)^B\subset \*c$, it is immediately clear that $j'$ is a $Z(B)\hbox{-}(\*c)^B$ bilinear map. The last assertion is proved as the second part of \leref{imj}.
\end{proof}

\begin{theorem}
\thlabel{denseB}
Let $\cc$ be an $A$-coring which is locally projective as a left $A$ module. 
Let $B\to A$ be a ring morphism and $I$ any index set. Consider $j\in{_A\Hom_{\*c}}(\cc,(\*c)^{I})$, denote $\tildej$ for the corresponding element in ${_{\cc^*}\Hom_A}(\cc^{(I)},\cc^*)$ 
and denote $\tildej':(\cc^{(I)})^B\to(\cc^*)^B$ for the restriction of $\tildej$.
Then the following statements are equivalent
\begin{enumerate}[(i)]
\item for all $c_1,\ldots,c_n\in\cc$ and $f\in(\cc^*)^B$, there exists an element $g\in \im \tildej'$ such that $g(c_i)=f(c_i)$ (i.e.\ $\im\tildej'$ is dense in the finite topology on $(\cc^*)^B$);
\item for all $c_1,\ldots,c_n\in\cc$, there exists an element $e\in \im \tildej'$ such that $e(c_i)=\varepsilon(c_i)$;
\item for all $c_1,\ldots,c_n\in\cc$ and $f\in(\cc^*)^B$, there exists an element $g\in \im \tildej'$  such that $f\cdot c_i= g\cdot c_i$ (i.e.\ $\im \tildej'$ is dense in the $\cc$-adic topology on $(\cc^*)^B$);
\item for all $c_1,\ldots,c_n\in\cc$, there exists an element $e\in \im \tildej'$ such that $e\cdot c_i= c_i$;
\item there exist $B$-linear local right inverses for $\tildej$, i.e.\ for all $c_1,\ldots,c_n\in\cc$, there exists a $\barj\in{_{\cc^*}\Hom_B}(\cc^*,(\cc)^{(I)})$ such that $\tildej(\barj(f))(c_i)=f(c_i)$ for all $f\in\cc^*$;
\item for all $c_1,\ldots,c_n\in\cc$, there exists an element $\bar{z}\in (\cc^B)^{(I)}$ such that $c_i=\bar{z}\cdot j(c_i)=\sum_\ell z_\ell j_\ell(c_i)$ for all $i=1,\ldots, n$.
\end{enumerate}
Moreover, if any of these conditions are satisfied, then
\begin{enumerate}[(a)]
\item there exist $B$-linear local left inverses for $j$, i.e.\ for all $c_1,\ldots,c_n\in\cc$, there exists a $\barj\in{_B\Hom_{\*c}((\*c)^{I},\cc)}$ such that $c_i=\barj(j(c_i))$ 
\item $j$ is injective; 
\item $\cc$ is $\im \tildej$-locally projective as a right $A$-module.
\end{enumerate}
\end{theorem}

\begin{proof}
$\ul{(i)\Rightarrow (ii)}$ Trivial.\\
$\ul{(ii)\Rightarrow (i)}$ By \leref{j'} we know that $\im \tildej'$ is a left ideal in $(\cc^*)^B$. The statement follows now immediately.\\
$\ul{(i)\Rightarrow (iii)}$ Take $c_i$ and $f$ as in statement $(iii)$. Then $f\cdot c_i=f(c_{i(1)})c_{i(2)}$. By statement $(i)$ we find an element $g\in\im\tildej'$ such that $g(c_{i(1)})=f(c_{(i(1))})$ for all $i$. Consequently, $g\cdot c_i=f\cdot c_i$.\\
$\ul{(iii)\Rightarrow (iv)}$ Follows again from the fact that $\im \tildej'$ is a left ideal in $(\cc^*)^B$.\\
$\ul{(iv)\Rightarrow (ii)}$ Take $c_i$ as in statement $(ii)$. From $(iv)$ we know that we can find an $e\in\im\tildej'$ such that $e\cdot c_i=e(c_{i(2)})c_{i(2)}=c_i$. Apply $\varepsilon$ to this last equation, then we find $e(c_i)=\varepsilon(c_i)$.\\
$\ul{(iv)\Rightarrow (vi)}$ Consider $c_1,\ldots,c_n\in\cc$. Then we know from $(iv)$ that there exists an element $e\in\im\tildej'$ such that $e\cdot c_i=c_i$. We can write $e=\tildej'(\bar{z})=\sum\tildej'_\ell(z_\ell)$ for some $\bar{z}=(z_\ell)\in (\cc^B)^{(I)}$.
We will show that this $\bar{z}$ is the needed one.
Recall from \prref{vcongv3} that $j_\ell$ is a right $\cc$-colinear morphism from $\cc$ to $\Rat_{\*c}(\*c)$. We find
\begin{eqnarray*}
\bar{z}\cdot j(c_i)&=&
\sum_\ell z_\ell\cdot j_\ell(c_i)=\sum_\ell {z_\ell}_{(1)}j_\ell(c_i)({z_\ell}_{(2)})\\
&=&\sum_\ell j_\ell(c_i)_{[0]}(z_\ell)j_\ell(c_i)_{[1]}
=\sum_\ell j_\ell({c_i}_{(1)})(z_\ell){c_i}_{(2)}\\
&=&\sum_\ell \tildej_\ell (z_\ell)({c_i}_{(1)}){c_i}_{(2)}
=e({c_i}_{(1)}){c_i}_{(2)}=c_i,
\end{eqnarray*}
where we used \leref{formulecf} in the third equation.\\
$\ul{(vi)\Rightarrow (v)}$ 
Take $c_i\in \cc$ as in the statement of $(v)$. Choose representatives $c_{j_i}, c_{k_i}\in \cc$ such that $\Delta(c_i)=\sum_ic_{k_i}\ot_A c_{j_i}$ for all $i$. By $(vi)$ we can find a 
$\bar{z}=(z_\ell)\in (\cc^B)^{(I)}$
such that for all $c_{k_i}$,
\begin{equation}\eqlabel{barj}
c_{k_i}=\bar{z}\cdot (j(c_{k_i})). 
\end{equation}
Now by \leref{modiso's} we can associate to $\bar{z}$ an element $\barj\in{_{\cc^*}\Hom_B}(\cc^*,\cc^{(I)})$, defined as $\barj(f)=f\cdot \bar{z}$ for all $f\in\cc^*$.
We find 
\begin{eqnarray*}
\tildej(\barj(f))(c_i)&=&\tildej(f\cdot \bar{z})(c_i)=(f*\tildej(\bar{z}))(c_i)
=\sum_\ell(f*\tildej_\ell(z_\ell))(c_i)\\
&=&\sum_{\ell,i} f(\tildej_\ell(z_\ell)(c_{k_i})c_{j_i})=\sum_{\ell,i}f(j_\ell(c_{k_i})(z_\ell)c_{j_i})\\
&=&\sum_if(\varepsilon(c_{k_i})c_{j_i})=f(c_i).
\end{eqnarray*}
Where the one but last equation follows by applying $\varepsilon$ on \equref{barj}.\\
$\ul{(v)\Rightarrow (i)}$ For every $f\in(\cc^*)^B$, we have $\barj(f)\in(\cc^B)^{(I)}$. Consequently we can choose $g=\tildej'(\barj(f))$.

Suppose now that the conditions $(i)-(vi)$ are satisfied. 
$\ul{(vi)\Rightarrow (a)}$ 
Follows immediately from \leref{modiso's}.
To prove $(b)$, suppose $j(c)=0$ for some $c\in\cc$, then by statement $(a)$ we can find $\barj$ such that $c=\barj(j(c))=\barj(0)=0$, so $j$ is injective.
Finally, we find by $(vi)$ on every set $c_1,\ldots,c_n\in\cc$ 
an element $(z_\ell)\in(\cc^B)^{(I)}$, such that we can compute 
$$c_i=\sum_\ell z_\ell\cdot j_\ell(c_i)
=\sum_\ell {z_\ell}_{(1)} j_\ell(c_i)({z_\ell}_{(2)})
=\sum_\ell {z_\ell}_{(1)}\tildej_\ell({z_\ell}_{(2)})(c_i).$$ 
This means that $\{{z_\ell}_{(1)},\tildej_\ell({z_\ell}_{(2)})\}$ is a local dual basis for $c_i$, so $\cc$ is locally projective as a right $A$-module.
\end{proof}

\begin{remark}\relabel{denseB}
It follows immediately from the proof that, even if $\cc$ is not necessary locally projective as a left $A$-module, the first four statements of \thref{denseB} remain equivalent statements if we replace $\im\tildej'$ by any left ideal in $(\cc^*)^B$.
\end{remark}

\begin{Definitions}
If $\cc$ is an $A$-coring that is locally projective as a left $A$-module and such that the equivalent conditions $(i)$-$(vi)$ of \thref{denseB} are satisfied we call $\cc$ a \emph{left $B$-locally quasi-Frobenius} coring. 

If $\cc$ is left $A$-locally quasi-Frobenius, we will just say that $\cc$ is \emph{left locally quasi-Frobenius}.

If $\cc$ is a $B$-locally quasi-Frobenius coring such that the index-set $I$ of \thref{denseB} can be chosen to have only 1 element, then we say that $\cc$ is \emph{left $B$-locally Frobenius}.
\end{Definitions}

\begin{corollary}\colabel{imjring}
Let $\cc$ be an $A$-coring that is left $B$-locally quasi-Frobenius with Frobenius morphism $j:\cc\to (\*c)^I$, denote as in \thref{denseB} the corresponding morphism $\tildej:\cc^{(I)}\to \cc^*$ with restricted morphism $\tildej':(\cc^{(I)})^B\to (\cc^*)^B$.
Then the following statements hold.
\begin{enumerate}[(i)]
\item $\im \tildej$ is a $B$-ring with left local units, where $\tildej:\cc^{(I)}\to \cc^*$. Moreover, $\im\tildej$ acts with local units on every left $\cc$-comodule;
\item $\im \tildej'$ is a ring with left local units and $\im\tildej'$ acts with local units on every left $\cc$-comodule.
\end{enumerate}
 
\end{corollary}

\begin{proof}
\ul{(i)}. Clearly $\im\tildej$ is a $B$-ring. Since $\im\tildej'\subset\im\tildej$, the remaining part of the statement follows by part (ii).\\
\ul{(ii)}.
Let $\tildej(c_i)$ be any element of $\im\tildej$ where $(c_i)\in \cc^{(I)}$. Then denote by $\bar{z}=(z_\ell)\in (\cc^B)^{(I)}$ the element satisfying condition $(vi)$ of \thref{denseB}. Write $e=\tildej'(z_\ell)=\tildej(z_\ell)\in\im\tildej'$, we claim that $e$ is a left local unit for $\tildej(c_i)$. Indeed,
\begin{eqnarray*}
\tildej(z_\ell)*\tildej(c_i)&=&\tildej\big(\tildej(z_\ell)\cdot c_i\big) 
= \tildej(\sum_\ell \tildej_\ell(z_\ell)\cdot c_i)\\
&=& \tildej(\sum_\ell \tildej_\ell(z_\ell)(c_{i(1)})c_{i(2)}) 
= \tildej(\sum_\ell j_\ell(c_{i(1)})(z_\ell)c_{i(2)}) \\
&=& \tildej(\sum_\ell j_\ell(c_i)_{[0]}(z_\ell) j_\ell(c_i)_{[1]})
= \tildej(\sum_\ell z_{\ell(1)}j_\ell(c_i)(z_{\ell(2)})) \\
&=& \tildej(\bar{z} \cdot j(c_i)) = \tildej(c_i).
\end{eqnarray*}
Here we used the left $\cc^*$-linearity of $\tildej$ in the first equality, \leref{formulecf} in the fifth equality and part $(vi)$ of \thref{denseB} in the last equality.

Let $M$ be any left $\cc$-comodule. The action of $f\in\im\tildej'$ on $m\in M$ is given by 
$f\cdot m=f(m_{[-1]})m_{[0]}$. That there exists local units for this action is a direct consequence of \thref{denseB}, part $(ii)$.
\end{proof}

\begin{theorem}\thlabel{QcoFroblocQFrob}
Let $\cc$ be an $A$-coring 
Then the following statements hold.
\begin{enumerate}
\item If $j\in{_A\Hom_{\*c}}(\cc,(\*c)^{I})$ is injective then the restriction $j':\cc^A\to ((\*c)^A)^{I}$ is also injective.
\item Consider ring morphisms $B\to B'\to A$. If $\cc$ left $B'$-locally quasi-Frobenius then $\cc$ is left $B$-locally quasi-Frobenius.
\item If $A$ is a PF-ring and $\cc$ is $\*c^*$-locally projective then the following statements are equivalent
\begin{enumerate}[(i)]
\item $\cc$ is left locally quasi-Frobenius;
\item $\cc$ is left $B$-locally quasi-Frobenius for an arbitrary ring morphism $B\to A$;
\item 
$\cc$ is left quasi-co-Frobenius;
\item The restriction $j':\cc^A\to((\*c)^A)^{I}$ of the Frobenius map $j$ is injective
(i.e.\ $\cc^A$ is a torsionless right $(\*c)^A$-module).
\end{enumerate}
\end{enumerate}
\end{theorem}

\begin{proof}
$\ul{(1)}$ Trivial.\\
$\ul{(2)}$ If $\cc$ is left $B'$-locally Frobenius, then the $B'$-linear local left inverses for $j$ from the equivalent condition $(vi)$ of \thref{denseB} are clearly also $B$-linear local left inverses.\\
$\ul{(3)}$ From part $(1)$ and $(2)$ we know already $\ul{(iii)\Rightarrow(iv)}$ and $\ul{(i)\Rightarrow(ii)}$. 
From \thref{denseB} we know that $\ul{(ii)\Rightarrow(iii)}$. So we only have to prove $\ul{(iv)\Rightarrow(i)}$. Let us denote $\tildej=\xi(j)\in {_{\cc^*}\Hom_A}(\cc^{(I)},\cc^*)$ by the isomorphism of \leref{modiso's}(i). We will show that $\im \tildej'$ is dense in the finite topology on $\*c^*$, which is equivalent condition $(i)$ of \thref{denseB} applied to the situation $B=A$. Since $\cc$ is $\*c^*$-locally projective, the canonical map $\cc\to (\*c^*)^*$ is injective. Moreover $A$ is a PF-ring, so a subset $P\subset\*c^*$ is dense in the finite topology if and only if the orthogonal complement $P^\perp$ of $P$ is trivial (see \cite[Theorem 1.8]{Abu:lintop}). 
Take any $c\in(\im\tildej')^\perp$. Then $\tildej'(d_\ell)(c)=j'(c)(d_\ell)=0$ for all $d_\ell\in\cc^{(I)}$. This implies $c\in\ker j'$. By the injectivity of $j'$ we find $c=0$, so $(\im\tildej')^\perp=0$ and $\im\tildej'$ is dense in the finite topology on $\*c^*$.
\end{proof}

\begin{corollary}\colabel{coFroblocFrob}
Let $\cc$ be an $A$-coring 
Then the following statements hold.
\begin{enumerate}
\item If $j\in{_A\Hom_{\*c}}(\cc,\*c)$ is injective then the restriction $j':\cc^A\to\*c$ is also injective.
\item Consider ring morphisms $B\to B'\to A$. If $\cc$ left $B$-locally Frobenius then $\cc$ is left $B'$-locally Frobenius.
\item If $A$ is a PF-ring and $\cc$ is $\*c^*$-locally projective then the following statements are equivalent
\begin{enumerate}[(i)]
\item $\cc$ is left co-Frobenius;
\item $\cc$ is left locally Frobenius;
\item $\cc$ is left $B$-locally Frobenius for an arbitrary ring morphism $B\to A$;
\item The restriction $j':\cc^A\to\*c$ of the Frobenius map $j$ is injective.
\end{enumerate}
\end{enumerate}
\end{corollary}

\begin{proof}
This is proved in the same way as \thref{QcoFroblocQFrob}.
\end{proof}

\begin{proposition}\prlabel{semiperfect}
Let $\cc$ be an $A$-coring and $B\to A$ any ring morphism. 
If $\cc$ is left $B$-locally quasi-Frobenius, then ${_{\cc^*}\Rat}$ is an exact functor. In particular, if $A$ is a QF-ring, then then $\cc$ is a left semiperfect coring.
\end{proposition}

\begin{proof}
By part $(2)$ of \coref{coFroblocFrob}, we know that $\cc$ is also $k$-locally quasi-Frobenius. This implies by \thref{denseB} that $\im \tildej$ is dense in $\cc^*$. Also by \thref{denseB}, we know that $\cc$ is locally projective as a right $A$-module, so \leref{modiso's}(i) implies that $\im \tildej$ is contained in ${_{\cc^*}\Rat}(\cc^*)$. We can conclude that ${_{\cc^*}\Rat}(\cc^*)$ itself is dense in $\cc^*$. By \cite[Proposition 2.6]{CI:semiperfect} the density of ${_{\cc^*}\Rat}(\cc^*)$ is equivalent to the exactness of ${_{\cc^*}\Rat}$. Moreover, if $A$ is a QF-ring, this condition is again equivalent to $\cc$ being a left semiperfect coring (see \cite[Theorem 4.3]{CI:semiperfect} or \cite[Theorem 3.8]{EGT:duality})
\end{proof}

\subsection{Characterization of Frobenius corings}\selabel{frobenius}

Considering the objects $\cc$ and $\*c$ in the category ${_A\Mm_{\*c}}$ and the objects $\cc$ and $\cc^*$ in the category ${_{\cc^*}\Mm_A}$, we obtain as in \seref{intrMorita} the following Morita contexts: 
\begin{eqnarray*}
\NN(\cc,\*c)&=&({_A\End_{\*c}}(\cc),{_A\End_{\*c}}(\*c),{_A\Hom_{\*c}}(\*c,\cc),{_A\Hom_{\*c}}(\cc,\*c),\circ,\bul);\\
\NN(\cc,\cc^*)&=&({_{\cc^*}\End_{A}}(\cc),{_{\cc^*}\End_{A}}(\cc^*),{_{\cc^*}\Hom_{A}}(\cc^*,\cc),{_{\cc^*}\Hom_{A}}(\cc,\cc^*),\bar{\circ},\bar{\bul}).
\end{eqnarray*}
If we consider the contravariant functors $\Jj$ and $\Kk,$ from \equref{JK} and the covariant functors $\widetilde{\Jj}$ and $\widetilde{\Kk}$ from \equref{JKtilde}, then we can construct another two Morita contexts
\begin{eqnarray*}
\YY(\Kk,\Jj)&=&(\Nat(\Kk,\Kk),\Nat(\Jj,\Jj),\Nat(\Jj,\Kk),\Nat(\Kk,\Jj),\newtriangleup,\newblacktriangleup);\\
\YY(\widetilde{\Kk},\widetilde{\Jj})&=&(\Nat(\widetilde{\Kk},\widetilde{\Kk}),\Nat(\widetilde{\Jj},\widetilde{\Jj}),\Nat(\widetilde{\Jj},\widetilde{\Kk}),\Nat(\widetilde{\Kk},\widetilde{\Jj}),
\widetilde{\newtriangleup},\widetilde{\newblacktriangleup}).
\end{eqnarray*}
Consider the functors $\Ff$, $\Gg$ and $\Hh$ as in \seref{inductionadjoint}. We can construct the Morita context that connects the functors $\Ff$ and $\Hh$ in the category of functors from $\Mm^\cc$ to $\Mm_A$ and all natural transformations between them.
$$\NN(\Ff,\Hh)=(\Nat(\Ff,\Ff),\Nat(\Hh,\Hh),\Nat(\Hh,\Ff),\Nat(\Ff,\Hh),\newsquare,\newblacksquare).$$
Although the functors $\Ff$ and $\Gg$ are not contained in the same category, we can apply the results of \seref{locadj} to obtain a Morita context \equref{MorContextFunctors} connecting the functors $\Ff$ and $\Gg$. 
$$\MM(\Ff,\Gg)=(\Nat(\Gg,\Gg),\Nat(\Ff,\Ff)^{\rm op},\Nat(\Mm_A,\Ff\Gg),\Nat(\Gg\Ff,\Mm^\cc),\newdiamond,\newblackdiamond).$$
Similarly, we find a Morita context connecting the functors $\Gg$ and $\Hh$,
$$\MM(\Gg,\Hh)=(\Nat(\Hh,\Hh),\Nat(\Gg,\Gg)^{\rm op}, \Nat(\Mm^\cc,\Gg\Hh),\Nat(\Hh\Gg,\Mm_A),\newtriangle,\newblacktriangle).$$

\begin{theorem}\thlabel{contextsFrob}\thlabel{Moritamorph}
Let $\cc$ be an $A$-coring. With notation as above, we have the following diagram of morphisms of Morita contexts.
\[
\xymatrix{
& \NN^{\rm top}(\Ff,\Hh) \ar[dr]^{{\nnn}} \ar[dl]_{\bar\nnn} \\
\MM(\Ff,\Gg) \ar[d]^\fff \ar[rr]_-{\mmm}  && \MM^{\rm op}(\Gg,\Hh) \\
\NN^{\rm top}(\cc,\cc^*) \ar@<.5ex>[rr]^{{\aaa}} \ar[d]_\bbb  && \NN(\cc,\*c) \ar[d]^{\bar{\bbb}}\\
\YY^{\rm top}(\Kk,\Jj) \ar[rr]^-{\bar\aaa} 
&& \YY(\widetilde{\Kk},\widetilde{\Jj})
}
\]
Here the upper script `op' indicates the opposite Morita context and `t' denotes the twisted Morita context (see \seref{intrMorita}).
For an arbitrary coring $\cc$, the morphisms $\aaa, \bar\aaa, \bbb, \bar\bbb, \mmm, \nnn, \bar\nnn$ are isomorphism of Morita contexts. 
If ${\cc}$ is locally projective as a left $A$-module, then $\fff$ becomes an isomorphism of Morita contexts as well.
\end{theorem}

\begin{proof}
The algebra isomorphisms for $\mmm, \nnn$ and $\bar\nnn$ follow immediately from \prref{XY} and \reref{XY}. The maps that describe the isomorphisms for the connecting bimodules are given in equations \equref{NatV} and \equref{NatW}.\\
The algebra isomorphisms for $\aaa, \bar\aaa, \bbb$ and $\bar\bbb$ follow from \prref{vw6} in combination with \reref{XY}. The isomorphisms for the connecting bimodules of $\aaa$ are given in \prref{V3V'3} and \prref{wcongw3}, for $\bar\aaa, \bbb$ and $\bar\bbb$ they follow from \prref{vw6}.\\
The first algebra morphism of $\fff$ is constructed as follows. We know by \reref{XY} that $\Nat(\Gg,\Gg)\cong{^\cc\End_A(\cc)}^{\rm op}$. Hence we have an algebra map 
$$\fff_1:\Nat(\Gg,\Gg)\cong{^\cc\End_A(\cc)}^{\rm op}\to {_{\cc^*}\End_A(\cc)}^{\rm op}.$$
The algebra map $\fff_2:\Nat(\Ff,\Ff)^{\rm op}\to {_{\cc^*}\End_A}(\cc^*)$ is given explicitly in \reref{XY}, and the bimodule maps $\fff_3$ and $\fff_4$ follow from \prref{wcongw3} and \prref{vcongv3} respectively. Moreover, we $\fff_2$ and $\fff_3$ are always bijective and when $\cc$ is flat as left $A$-module, then $\fff_1$ is an isomorphism by a rationality argument and $\fff_4$ is an isomorphism by \prref{vcongv3}.\\
We leave it to the reader to verify that all given (iso)morphisms of algebras and bimodules do indeed form Morita morphisms and that the stated diagrams of Morita morphisms commute.
\end{proof}

\begin{corollary}\colabel{fingen}
Let $\cc$ be an $A$-coring. Then there exists a split epimorphism $j\in{_{\cc^*}\Hom_A}(\cc,\cc^*)$ if and only if there exists a split monomorphism $\tildej\in{_A\Hom_{\*c}}(\cc,\*c)$. If any of these equivalent conditions holds then $\cc$ finitely generated and projective as a right $A$-module.
\end{corollary}

\begin{proof}
We will prove a more general version of this corollary in \coref{fingenQuasi}
\end{proof}

As a corollary we obtain the well-known characterization of Frobenius corings in terms of Frobenius functors.

\begin{corollary}[\textbf{characterization of Frobenius corings}]\colabel{Frobenius}
Let $\cc$ be an $A$-coring, then the following statements are equivalent;
\begin{enumerate}[(i)]
\item $\cc\cong \*c$ in ${_A\Mm_{\*c}}$ (i.e.\ $\cc$ is a Frobenius coring);
\item $\cc\cong \cc^*$ in ${_{\cc^*}\Mm_A}$;
\item the functors $\Hh$ and $\Ff$ are naturally isomorphic;
\item $(\Gg,\Ff)$ is a pair of adjoint functors, and therefore $(\Ff,\Gg)$ is a Frobenius pair;
\item $(\Hh,\Gg)$ is a pair of adjoint functors, and therefore $(\Gg,\Hh)$ is a Frobenius pair;
\item the functors $\Jj$ and $\Kk$ are naturally isomorphic;
\item the functors $\widetilde{\Jj}$ and $\widetilde{\Kk}$ are naturally isomorphic;
\item ${_A\cc}$ is finitely generated and projective and the functors $\Jj'$ and $\Kk'$ are isomorphic;
\item left hand versions of $(iii)-(viii)$, replacing $\Mm^\cc$ by ${^\cc\Mm}$ and $\Mm_A$ by ${_A\Mm}$. 
\end{enumerate}
\end{corollary}

\begin{proof}
The first statement is true if and only if there exists a pair of invertible elements in the Morita context $\NN(\cc,\*c)$. From (the left hand version of) \coref{fingen} we know that the isomorphism $\cc\cong\cc^*$ implies that $\cc$ is finitely generated and projective as a left and right $A$-module. The equivalence of $(i)-(vii)$ follows now immediately from the (anti-)isomorphisms of Morita contexts from \thref{contextsFrob}. Note that $(\Ff,\Gg)$ is always a pair of adjoint functors and therefore the adjointness of $(\Gg,\Ff)$ means exactly that $(\Ff,\Gg)$ is a Frobenius pair. The same reasoning holds for the pair $(\Gg,\Hh)$.
Since for a coring that is finitely generated and projective as a left $A$-module the categories $\Mm^\cc$ and $\Mm_{\*c}$ are isomorphic, we obtain that the functor $\Jj'$, respectively $\Kk'$, is isomorphic with $\Jj$, respectively $\Kk$.
The equivalence of $(i)$ and $(ii)$ imply the equivalence with the left hand version of the other statements. 
\end{proof}

\begin{remark}\relabel{symmetry}
\begin{enumerate}
\item
Note that the left-right symmetry of the 
notion of a Frobenius extension (or a Frobenius coring), is by the previous Corollary a consequence of the isomorphism of Morita contexts between $\NN(\cc,\*c)$ and $\NN^{\rm top}(\cc,\cc^*)$. We will see that this isomorphism is missing in \thref{Scontexts} if we study the quasi-co-Frobenius property in \seref{QcF}. This (partially) explains why the notion of a quasi-co-Frobenius coring is not left-right symmetric.
\item 
Considering the functors $\Jj'$ and $\Kk'$ (see \equref{JK}), we can construct another Morita context $\YY(\Kk',\Jj')$. A pair of invertible elements in this context describes when 
$\cc\cong {_\Rr\Rat(\cc^*)}$. This can be in particular of interest when $\cc=C$ is a coalgebra over a field, since in that case we know from \cite{mio} $C$ is at the same time left and right co-Frobenius if and only if $C\cong \Rat(C^*)$.
\end{enumerate}
\end{remark}

\subsection{Quasi-co-Frobenius corings and related functors}\selabel{QcF}

Let $I$ be any index set and consider the objects $\cc$ and $(\*c)^I$ in the category ${_A\Mm_{\*c}}$ and the objects $(\cc)^{(I)}$ and $\cc^*$ in the category ${_{\cc^*}\Mm_A}$, we obtain in this way the Morita contexts 
\begin{eqnarray*}
\NN((\*c)^I,\cc)&=&({_A\End_{\*c}}((\*c)^I),{_A\End_{\*c}}(\cc),{_A\Hom_{\*c}}(\cc,(\*c)^I),{_A\Hom_{\*c}}((\*c)^I,\cc),\circ,\bul);\\
\NN(\cc^{(I)},\cc^*)&=&({_{\cc^*}\End_{A}}(\cc^{(I)}),{_{\cc^*}\End_{A}}(\cc^*),{_{\cc^*}\Hom_{A}}(\cc^*,\cc^{(I)}),{_{\cc^*}\Hom_{A}}(\cc^{(I)},\cc^*),\bar{\circ},\bar{\bul}).
\end{eqnarray*}
Consider again the functors $\Jj$ and $\widetilde{\Jj}$ from \equref{JK} and \equref{JKtilde} and the functors $\Kk^s$ and $\widetilde{\Kk}^s$ from \equref{KkS}. We can construct the following Morita contexts.
\begin{eqnarray*}
\YY(\Kk^s,\Jj)&=&(\Nat(\Kk^s,\Kk^s),\Nat(\Jj,\Jj),\Nat(\Jj,\Kk^s),\Nat(\Kk^s,\Jj),\newtriangleup,\newblacktriangleup);\\
\YY(\widetilde{\Kk}^s,\widetilde{\Jj})&=&(\Nat(\widetilde{\Kk}^s,\widetilde{\Kk}^s),\Nat(\widetilde{\Jj},\widetilde{\Jj}),\Nat(\widetilde{\Jj},\widetilde{\Kk}^s),\Nat(\widetilde{\Kk}^s,\widetilde{\Jj}),\widetilde{\newtriangleup},\widetilde{\newblacktriangleup}).
\end{eqnarray*}
Dually, we can consider functors $\Ii={_A\Hom_{\*c}}(-,\*c), \Ll^s={_A\Hom_{\*c}}(-,\cc^{(I)}):{_A\Mm_{\*c}}\to{\Mm_{\*c^*}}$ and 
$\widetilde\Ii={_A\Hom_{\*c}}(\*c,-), \widetilde\Ll^s={_A\Hom_{\*c}}(\cc^{(I)},-):{_A\Mm_{\*c}}\to{_{\*c^*}\Mm}$. Out of these functors we construct Morita contexts $\YY(\Ll^s,\Ii)$ and $\YY(\widetilde\Ll^s,\widetilde\Ii)$.\\
Consider the functors $\Ff,\Gg,\Hh$ and $\Ss$ from \seref{induction}. 
We immediately obtain the following Morita context.
$$\NN(\Ss\Ff,\Hh)=(\Nat(\Ss\Ff,\Ss\Ff),\Nat(\Hh,\Hh),\Nat(\Hh,\Ss\Ff),\Nat(\Ss\Ff,\Hh),\newsquare,\newblacksquare).$$
Applying the techniques of \seref{locadj}, we find a Morita context of type 
\equref{MorContextFunctors} connecting $\Ff$ and $\Gg\Ss$
and a context connecting $\Ss\Ff$ and $\Gg$.
\begin{eqnarray*}
\MM(\Ff,\Gg\Ss)&=&
(\Nat(\Gg\Ss,\Gg\Ss),\Nat(\Ff,\Ff)^{\rm op},\Nat(\Mm_A,\Ff\Gg\Ss),\Nat(\Gg\Ss\Ff,\Mm^\cc),\newdiamond,\newblackdiamond);
\\
\MM(\Ss\Ff,\Gg)&=&
(\Nat(\Gg,\Gg),(\Nat(\Ss\Ff,\Ss\Ff)^{\rm op},\Nat(\Mm_A,\Ss\Ff\Gg),\Nat(\Gg\Ss\Ff,\Mm^\cc),\bar\newdiamond,\bar\newblackdiamond).
\end{eqnarray*}
Let us give the explicit form of the connecting maps. Denote $\alpha\in\Nat(\Gg\Ss\Ff,\Mm^\cc)$, 
$\beta\in\Nat(\Mm_A,\Ff\Gg\Ss), \beta'\in\Nat(\Mm_A,\Ss\Ff\Gg)$ and $\gamma\in\Nat(\Gg\Ss\Ff,\Mm^\cc)$, $N\in\Mm_A$ and $M\in\Mm^\cc$ then
\begin{eqnarray*}
(\beta{\newdiamondd}\gamma)_N = \gamma_{\Gg\Ss N}\circ\Gg\Ss\beta_N, &&
(\gamma{\newblackdiamondd}\beta)_M = \Ff\gamma_M\circ \beta_{\Ff M};\\
(\beta'\bar{\newdiamondd}\alpha)_N=\alpha_{\Gg N}\circ\Gg\beta'_N, && 
(\alpha\bar{\newblackdiamondd}\beta')_N=\Ss\Ff\alpha_M\circ\beta'_{\Ss\Ff M}. \\
\end{eqnarray*}
The left hand versions of the previous two contexts can be obtained by considering the functors ${^\cc\Ff}=\Ff', {^\cc\Hh}=\Hh':{^\cc\Mm}\to{_A\Mm}$, ${^\cc\Gg=\Gg'}:{_A\Mm}\to{^\cc\Mm}$ and $\Ss':{_A\Mm}\to{_A\Mm}$. This way, we obtain contexts $\NN(\Ss'\Ff',\Hh')$, $\MM(\Ff',\Gg'\Ss')$ and $\MM(\Ss'\Ff',\Gg')$.
 
\begin{theorem}\thlabel{Scontexts}
Let $\cc$ be an $A$-coring and keep the notation from above. 
\begin{enumerate}[(i)]
\item
There exist morphisms of Morita contexts as in the following diagram.
\[
\xymatrix{
&& \NN(\Ss'\Ff',\Hh') \ar[dl]_{\nnn} \ar[dr]^{\nnn'}\\
&\MM(\Ff,\Gg\Ss) \ar[rr]^{\mmm} \ar[d]^{\fff} && \MM^{\rm top}(\Ss'\Ff',\Gg')\\
& \NN(\cc^{(I)},\*c) \ar[dr]^{\bbb} \ar[dl]_{\bbb'}\\
\YY(\Ll^s,\Ii) \ar[rr]^{\ccc} && \YY(\widetilde\Ll^s,\widetilde\Ii)
}
\]
For an arbitrary coring $\cc$, the morphisms $\bbb, \bbb', \ccc$ and $\mmm, \nnn, \nnn'$ are isomorphisms of Morita contexts, if $\cc$ is locally projective as left $A$-module, then $\fff$ is an isomorphism of Morita contexts as well.
\item
There exist morphisms of Morita contexts as in the following diagram.
\[
\xymatrix{
&& \NN(\Ss\Ff,\Hh) \ar[dl]_{\bar\nnn} \ar[dr]^{\bar\nnn'}\\
& \MM(\Ff',\Gg'\Ss') \ar[rr]^{\bar\mmm} \ar[d]^{\bar\fff} && \MM^{\rm top}(\Ss\Ff,\Gg) \\
& \NN(\cc^{(I)},\cc^*) \ar[dl]_{\bar\bbb'} \ar[dr]^{\bar\bbb}\\
\YY(\Kk^s,\Jj) \ar[rr]^{\bar\ccc} && \YY^{\rm op}(\widetilde\Kk^s,\widetilde\Jj) 
}
\]
For an arbitrary coring $\cc$, the morphisms $\bar\bbb, \bar\bbb', \bar\ccc$ and $\bar\mmm, \bar\nnn, \bar\nnn'$ are isomorphisms of Morita contexts, if $\cc$ is locally projective as left $A$-module, then $\bar\fff$ is an isomorphism of Morita contexts as well.
\item
There exists a anti-morphism of Morita contexts
\[
\xymatrix{
\NN(\cc^{(I)},\cc^*) \ar[rr]^{\aaa} && \NN(\cc,(\*c)^I). 
}
\]
\end{enumerate}
\end{theorem}

\begin{proof}
\ul{(i)}.
The algebra isomorphisms for $\nnn$, $\nnn'$ and $\mmm$ follow immediately from \prref{XY}, \reref{XY} and \prref{SFSF}. The maps that describe the isomorphisms for the connecting bimodules are given in \prref{VWS}, together with their left-right dual versions.\\
The Morita morphisms $\bbb$ and $\bbb'$ are obtained by left-right duality out of $\bar\bbb$ and $\bar\bbb'$ of part (ii). To construct $\bar\bbb$ we can work as follows.
The algebra isomorphism $\bar\bbb_2$ follows from \prref{vw6} in combination with \reref{XY}. The algebra isomorphism $\bar\bbb_1$ is given in \prref{NatYonS}.
The isomorphisms for the connecting bimodules of $\bar\bbb_3$ and $\bar\bbb_4$ are given in \prref{NatYonS}. The morphism $\bar\bbb'$ is constructed in the same way. \\
The first algebra morphism of $\fff$ is constructed as follows. We know by \prref{SFSF} that $\Nat(\Gg\Ss,\Gg\Ss)\cong{^\cc\End_A(\cc^{(I)})}^{\rm op}$. Hence we have an algebra map 
$$\fff_1:\Nat(\Gg\Ss,\Gg\Ss)\cong{^\cc\End_A(\cc^{(I)})}^{\rm op}\to {_{\cc^*}\End_A(\cc^{(I)})}^{\rm op}.$$
The algebra map $\fff_2:\Nat(\Ff,\Ff)^{\rm op}\to {_{\cc^*}\End_A}(\cc^*)$ is given explicitly in \reref{XY}, and the bimodule maps $\fff_3$ and $\fff_4$ follow from \leref{modiso's} and \prref{SFG,GSF} respectively. Moreover, we $\fff_2$ and $\fff_3$ are always bijective and when $\cc$ is locally projective as left $A$-module, then $\fff_1$ is an isomorphism by a rationality argument and $\fff_4$ is an isomorphism by \leref{modiso's}.\\
We leave it to the reader to verify that all given (iso)morphisms of algebras and bimodules do indeed form Morita morphisms and that the stated diagrams of Morita morphisms commute.
\\
\ul{(ii)}. Follows by left-right duality.\\
\ul{(iii)}.
Consider the element $\bar{\varepsilon}\in \Hom(\cc^{(I)},A)$, defined by the following diagram,
\[
\xymatrix{
\cc^{(I)} \ar[rr]^{\bar{\varepsilon}} && A\\
\cc \ar[u]^{\iota_\ell} \ar[urr]_\varepsilon 
}
\]
Then we have map ${_{\cc^*}\End_A}(\cc^{(I)})\to{_A\Hom_A}(\cc^{(I)},A)$ defined by composing with $\bar{\varepsilon}$ on the left.
Combining \leref{sumproduct} and \leref{M^B} we find that 
$${_A\Hom_A}(\cc^{(I)},A)\cong ({_A\Hom_A}(\cc,A))^I={^*\cc^*}^I\cong ((\*c)^I)^A\cong {_A\Hom_{\*c}}((\*c)^I,(\*c)^I).$$ 
Composing these maps we obtain a linear map $\aaa_1:{_{\cc^*}\End_A}(\cc^{(I)})\to {_A\End_{\*c}}((\*c)^I)$ one can easily check that this is an anti-algebra morphism. 
The algebra map $\aaa_2:{_{\cc^*}\End_A}(\cc^*)\to{_A\End_{\*c}}(\cc)^{\rm op}$ is constructed as in the proof of \thref{Moritamorph} part (i).
From \leref{sumproduct} we obtain an isomorphism $\aaa_3:{_{\cc^*}\Hom_{A}}(\cc^{(I)},\cc^*)\to{_A\Hom_{\*c}}(\cc,(\*c)^I)$.
The last morphism $\aaa_4:{_{\cc^*}\Hom_{A}}(\cc^*,\cc^{(I)})\to {_A\Hom_{\*c}}((\*c)^I,\cc)$ is constructed as follows. Denote $f(\varepsilon)=(z_\ell)$ for any $f\in{_{\cc^*}\Hom_{A}}(\cc^*,\cc^{(I)})$. Then we define $\aaa_4(f)(f_\ell)=\sum_\ell z_\ell\cdot f_\ell$, for $(f_\ell)\in (\*c)^I$.
The reader can check that the four morphisms together make up an anti-morphism of Morita contexts.
\end{proof}

\begin{corollary}\colabel{fingenQuasi}
Let $\cc$ be an $A$-coring, then there exists a split epimorphism $j\in{_{\cc^*}\Hom_A}(\cc^{(I)},\cc^*)$ if and only if there exists a split monomorphism $\tildej\in{_A\Hom_{\*c}}(\cc,(\*c)^I)$, whose left inverse is induced by an element $(z_\ell)\in (\cc^A)^{(I)}$. If any of these equivalent conditions holds then $\cc$ finitely generated and projective as a right $A$-module.
\end{corollary}

\begin{proof}
Consider the anti-morphism of Morita contexts $\aaa$ of \thref{Scontexts}(iii). 
First note that the condition for the left inverse of $\tildej$ means exactly that it lies inside the image of $\aaa_4$.
Suppose $j$ has a right  inverse $\barj$. Consider the morphism of Morita contexts $\aaa$ from \thref{Scontexts}. Then we obtain that $\aaa_4(\barj)$ is a left inverse for $\aaa_3(j)$. For the converse, suppose that $\tildej$ has a left inverse of the form $\aaa_4(\barj)$. We know that $\aaa_3$ is an isomorphism, so we can write $\tildej=\aaa_3(j)$ for some morphism $j\in{_{\cc^*}\Hom_A}(\cc^{(I)},\cc^*)$. Then we find $\aaa_4(\barj)\circ\aaa_3(j)=\cc=\aaa_2(\cc^*)$. Since $\aaa$ is an anti-morphism of Morita contexts we find that $\barj$ is a right inverse for $j$. Finally, denote $(z_\ell)\in(\cc^A)^{(I)}$ for the representative of the left inverse of $\tildej$. Then we find for all $c\in\cc$,
$$c=\sum_\ell z_{\ell(1)}\tildej_\ell(c)(z_{\ell(2)})
=\sum_\ell z_{\ell(1)}j_\ell(z_{\ell(2)})(c),$$
i.e.\ $\{z_{\ell(1)},j_\ell(z_{\ell(2)})\}$ is a finite dual basis for $\cc$ as a right $A$-module.
\end{proof}

\begin{theorem} \thlabel{QcF}
Suppose that $\cc$ is an $A$-coring which is locally projective as a left $A$-module. Then
\begin{enumerate}[(i)]
\item The following statements are equivalent
\begin{enumerate}[(a)]
\item $\cc$ is left locally quasi-Frobenius; 
\item there exists a $\cc^*$-$A$ bilinear map $\tildej:\cc^{(I)}\to\cc^*$ 
such that 
$$T(\tildej)=\{\tildej\circ\psi ~|~ \psi\in{_{\cc^*}\Hom_A}(\cc^*,\cc^{(I)}) \}\subset 
{_{\cc^*}\End_A}(\cc^*)\cong (\*c^*)^{\rm op}$$ 
acts unital on all objects of the generating subcategory ${^\cc_{\rm fgp}\Mm}$ of ${^\cc\Mm}$;
\item there exists a $\cc^*$-$A$ bilinear map $\tildej:\cc^{(I)}\to\cc^*$ 
such that $T(\tildej)$ acts with right local units on $\cc$;
\item there exists a natural transformation $J':\Gg'\Ss'\Ff'\to {^\cc\Mm}$ such that 
$$R(J')=\{J'{\newblackdiamondd'} \beta' ~|~ \beta'\in \Nat({_A\Mm},\Ff'\Gg'\Ss') \}\subset \Nat(\Ff',\Ff')^{\rm op}$$ 
acts unital on the generating subcategory ${^\cc_{\rm fgp}\Mm}$ of ${^\cc\Mm}$;
\item $\Ff'$ is a right ${^\cc_{\rm fgp}\Mm}$-locally quasi-adjoint for $\Gg'$;
\item there exists a natural transformation $J:\Gg\Ss\Ff\to {\Mm^\cc}$ such that 
$$S(J)=\{\beta{\bar\newdiamondd} J ~|~ \beta\in \Nat(\Mm_A,\Ss\Ff\Gg) \}\subset \Nat(\Gg,\Gg)^{\rm op}$$ 
acts unital on the generating subcategory ${\Mm^\cc_{\rm fgp}}$ of ${\Mm^\cc}$;
\item $\Gg$ is a left ${\Mm^\cc_{\rm fgp}}$-locally quasi-adjoint pair for $\Ff$;
\item there exists a natural transformation $\alpha\in\Nat(\Kk^s,\Jj)$ such that $$\{\alpha\circ\beta~|~\beta \in\Nat(\Jj,\Kk^s)\}\subset\Nat(\Jj,\Jj)\cong {\*c^*}^{\rm op}$$ 
acts with right local units on $\cc$;
\item[(j)] there exists a natural transformation $\alpha\in\Nat(\widetilde{\Jj},\widetilde{\Kk}^s)$ such that 
$$\{\beta\circ\alpha~|~\beta \in\Nat(\widetilde{\Kk}^s,\widetilde{\Jj})\}\subset\Nat(\widetilde{\Jj},\widetilde{\Jj})\cong {\*c^*}$$ 
acts with left local units on $\cc$;
\end{enumerate}
\item dually, we can characterize right locally quasi-Frobenius corings; in particular $\cc$ is right locally quasi-Frobenius if and only if $\Gg'$ is a left $\Mm^\cc_{\rm fgp}$-locally quasi-adjoint for $\Ff'$ if and only if $\Ff$ is a right ${^\cc_{\rm fgp}\Mm}$-locally quasi-adjoint for $\Gg$;
\item $\cc$ is left and right locally quasi-Frobenius if and only if 
$(\Ff,\Gg)$ is a $\Mm^\cc_{\rm fgp}$-locally quasi-Frobenius pair of functors if and only if $(\Ff',\Gg')$ is a ${^\cc_{\rm fgp}\Mm}$-locally quasi-Frobenius pair of functors.
\end{enumerate}
\end{theorem}

\begin{proof}
$\ul{(i)}$. 
$\ul{(a)\Rightarrow (b)}$ Suppose that $\cc$ is left locally Frobenius. We know by \coref{imjring} that $\im\tildej'$, where $\tildej':(\cc^{(I)})^A\to({\cc^*})^A$, acts with (left) local units on the objects of ${^\cc\Mm}$, and therefore unital on the objects of ${^\cc_{\rm fgp}\Mm}$ (see \thref{locunit}). Since ${_{\cc^*}\Hom_A}(\cc^*,\cc^{(I)})\cong (\cc^{(I)})^A$ (\leref{modiso's}), we can identify $\im\tildej'$ with $T(\tildej)^{\rm op}$ and the statement follows.\\
$\ul{(b)\Rightarrow (c)}$. Follows by \thref{locunit}.\\
$\ul{(c)\Rightarrow (a)}$. Follows by \thref{denseB} and \coref{imjring} using the same interpretation of $T(\tildej)$ as in the proof of part $\ul{(a)\Rightarrow (b)}$.\\
$\ul{(b)\Leftrightarrow (d)}$. Since a left locally quasi-Frobenius coring is locally projective as left $A$-module, this is in fact an immediate consequence of the isomorphism of Morita contexts $\bar\fff$ of \thref{Scontexts}, part (ii). We give however a direct proof.\\
Condition $(d)$ means that for any $M\in{^\cc_{\rm fgp}\Mm}$ and any left $A$-module morphism $f:\Ff' M\to \Ff' M'$ with $M'\in{^\cc\Mm}$, we can find a $\beta\in\Nat({_A\Mm},\Ff'\Gg'\Ss')$ such that
$$(J'{\newblackdiamondd}' \beta)_{M'}=\Ff' J'_{M'}\circ\beta_{\Ff' M'} \circ f = f$$
i.e.\ the following diagram commutes
\begin{equation}\eqlabel{Sadj1}
\xymatrix{
& \Ff' M \ar[dr]^f \ar[dl]_f \\
\Ff' M' \ar[rr]^{(J'{\newblackdiamondd}' \beta)_{M'}} 
\ar[dr]_{\beta_{\Ff' M'}} && \Ff' M'\\
& \Ff'\Gg'\Ss'\Ff' M' \ar[ur]_{\Ff' J'_{M'}}
}
\end{equation}
where the commutativity of the lower triangle is nothing else than the definition of ${\newblackdiamond}'$.
Since $\cc$ is locally projective as a left $A$-module, we find by \thref{Scontexts} an isomorphism of Morita contexts $\bar\fff:\NN(\cc^{(I)},\cc^*)\to \MM(\Ff',\Gg'\Ss')$.
This implies that that the existence of a natural transformation $\beta$ is equivalent to the existence of a morphism $\psi\in{_{\cc^*}\Hom_A}(\cc^*,\cc^{(I)})$. We can translate diagram \equref{Sadj1} now into the following diagram
\[
\xymatrix{
& M \ar[dl]_f \ar[dr]^f\\
M' \ar[dr]_{\psi_1} && M' \\
& \cc^{(I)}\ot_A M' \ar[ur]_{\psi_2}
}
\]
where $\psi_1$ and $\psi_2$ are given by
\begin{eqnarray}
\psi_1(m) &=&  \psi(\varepsilon)\ot_A m\\
\psi_2((c_\ell)\ot_A m ) &=&  \tildej(c_\ell)\cdot m,
\end{eqnarray}
where $m\in M'$ and $c_\ell\in \cc^{(I)}$. So the above diagram commutes if and only if 
$m=\tildej(\psi(\varepsilon))\cdot m$ for all $m\in \im f$, i.e.\ if and only if we can find a local unit for all elements of $\im f$ and this local unit has to be of the form $\tildej\circ\psi(\varepsilon)$. Note that this local unit is exactly an element of $\im\tildej'\cong T(\tildej)^{\rm op}$. If condition $(b)$ holds, then we know that there exists such a unit for all left $\cc$-comodules that are finitely generated and projective as a left $A$-module, so in particular we find a local unit for $\im f$, and thus condition $(d)$ holds as well. Conversely, if condition $(d)$ holds, than we find as above a local unit in $\im\tildej'\cong T(\tildej)^{\rm op}$ for all modules of the form $\im f$. Taking $M'=M$ and $f$ the identity map, we obtain a local unit for all $M\in {^\cc_{\rm fgp}\Mm}$, i.e.\ $(b)$ is satisfied as well.\\
$\ul{(d)\Leftrightarrow (e)}$. Follows directly from the definition.\\
$\ul{(c)\Leftrightarrow (f)}$. 
Condition $(f)$ means that for any $M\in\Mm^\cc_{\rm fgp}$ and $f:M\to \Gg(N)=N\ot_A\cc$ with $N\in\Mm_A$, there exists $\beta\in\Nat(\Mm_A,\Ss\Ff\Gg)$ such that 
$$(\beta\bar\newdiamondd J)_N\circ f=J_{\Gg N}\circ \Gg\beta_N\circ f=f,$$
or the following diagram commutes.
\[
\xymatrix{
& M \ar[dl]_f \ar[dr]^f\\
\Gg(N) \ar[rr]^-{(\beta\bar\newdiamondd J)_N} \ar[dr]_{\Gg\beta_N} && \Gg(N)\\
&\Gg\Ss\Ff\Gg(N) \ar[ur]_{J_{\Gg N}}
}
\]
Since $\cc$ is locally projective as a left $A$-module, we find by \thref{Scontexts} an isomorphism between the Morita contexts $\NN(\cc^{(I)},\cc^*)$ and $\MM^{\rm top}(\Ss\Ff,\Gg)$.
Thus, the existence of $\beta$ as above is equivalent to the existence of an $\cc^*$-$A$ bilinear map $\psi:\cc^*\to\cc^{(I)}$ such that the following diagram commutes
\begin{equation}\eqlabel{Sadj2}
\xymatrix{
& M \ar[dl]_f \ar[dr]^f\\
N\ot_A\cc \ar[dr]_{\psi_1} && N\ot_A\cc\\
& N\ot_A\cc^{(I)}\ot_A\cc \ar[ur]_{\psi_2}
}
\end{equation}
where $\psi_1$ and $\psi_2$ are given by
\begin{eqnarray}\eqlabel{psi12}
\psi_1(n\ot_A c)&=&n\ot_A \psi(\varepsilon)\ot_A c\\
\psi_2(n\ot (c_i)\ot c)&=&n\ot_A \tildej(c_i)\cdot c.
\end{eqnarray}
Here we denoted $n\in N$, $c\in\cc$ and $(c_i)\in\cc^{(I)}$. Then diagram \equref{Sadj2} commutes if and only if $\sum_i n_i\ot \tildej(\psi(\varepsilon))\cdot c_i =\sum_i n_i\ot_A c_i$ for all $\sum_i n_i\ot_A c_i \in \im f$. 
Suppose that condition $(f)$ holds and take any $c\in\cc$. Put $N=A$ and $M=cA$, the cyclic right $A$-module generated by $c$ and let $f:M\to A\ot_A\cc\cong \cc$ be the canonical injection.
Then by diagram \equref{Sadj2}, we obtain a left local unit $\tildej\circ\psi(\varepsilon)$ on $M$, i.e.\ we find left local unit in $T(\tildej)^{\rm op}$ for $c$. This shows that $(f)$ implies $(c)$.\\
Conversely, if condition $(c)$ is satisfied, then we know that we can find a left local unit in $T(\tildej)$ for any finite number of elements in $\cc$.
Take any $M\in\Mm^\cc_{\rm fgp}$ and $f:M\to N\ot_A\cc$. Then $\im f$ is also finitely generated. Take a finite number of generators for $\im f$ and denote representatives of them by
$n^i\ot_Ac^i$ (to reduce the number of indices, we omit a summation if we denote an element of $N\ot_A\cc$). 
By $(c)$ we know that we can find a left local unit $e=\tildej\circ\psi(\varepsilon)\in T(\tildej)^{\rm op}$ for the generators $c^i$, i.e.\ such that $e\cdot c^i=c^i$ for all $i$, for a particular choice of $\psi\in{_{\cc^*}\Hom_A(\cc^*,\cc^{(I)})}$. If we define $\psi_1$ and $\psi_2$ as in \equref{psi12}, we find that diagram \equref{Sadj2} commutes and we obtain indeed that $(c)$ implies $(f)$.\\
$\ul{(f)\Leftrightarrow (g)}$. Follows directly from the definition.\\
$\ul{(c)\Leftrightarrow (h) \Leftrightarrow (j)}$. Follows from the isomorphisms of Morita contexts $\bar\bbb$ and $\bar\bbb'$ of \thref{Scontexts}.\\
\ul{(ii)}. Follows by left-right duality.\\
\ul{(iii)}. Is a direct combination of the first two parts.
\end{proof}

\begin{remark}
If one takes the index-set $I$ to contain a single element in the previous Theorem, then we obtain a characterization of locally Frobenius corings (and consequently of co-Frobenius corings if the base ring is a PF-ring). In particular, we find that an $A$-coring $\cc$ is left locally Frobenius if and only if $\Ff'$ is a right ${^\cc_{\rm fgp}\Mm}$-locally adjoint for $\Gg'$ if and only if $\Gg$ is a left $\Mm^\cc_{\rm fgp}$-locally adjoint for $\Ff$. Moreover $\cc$ is at the same time left and right locally Frobenius if and only if $(\Ff,\Gg)$ is a $\Mm^\cc_{\rm fgp}$-locally Frobenius pair if and only if $(\Ff',\Gg')$ is a ${^\cc_{\rm fgp}\Mm}$-locally Frobenius pair.
\end{remark}

\end{document}